\documentclass[11pt]{article}
\input{glyphtounicode}
\usepackage[T1]{fontenc}
\usepackage[margin=1in]{geometry}
\DeclareUnicodeCharacter{2261}{\ensuremath{\equiv}}
\usepackage{lmodern}
\usepackage{graphicx}
\usepackage{xcolor}
\usepackage{multicol,multirow}
\usepackage{rotating}
\usepackage{longtable}
\usepackage[sort&compress,numbers]{natbib}
\usepackage[breaklinks,colorlinks,allcolors=blue]{hyperref}
\usepackage{breakurl}
\usepackage{amssymb,amsmath,amsthm,latexsym}
\usepackage{tikz}
\usetikzlibrary{arrows.meta,positioning,calc,fit,backgrounds,shapes.geometric}
\theoremstyle{definition}
\newtheorem{definition}{Definition}
\theoremstyle{plain}
\newtheorem{lemma}{Lemma}
\newtheorem{theorem}{Theorem}

\theoremstyle{definition}

\usepackage{coqdoc}
\usepackage{microtype}
\newcommand{\coq}[1]{\mbox{\normalfont\ttfamily #1}}
\newcommand{\coqcode}[2]{\coq{#2}}
\newcommand{\SameCenter}{\coqcode{same_center}{same\_center}}
\newcommand{\PointSetEq}{\coqcode{PointSet_eq}{PointSet\_eq}}
\newcommand{\PointSetIncl}{\coqcode{PointSet_incl}{PointSet\_incl}}
\newcommand{\PointSetStable}{\coqcode{PointSet_stable}{PointSet\_stable}}
\newcommand{\PointSetUnion}{\coqcode{PointSet_union}{PointSet\_union}}

\newcommand{\PointSetMeet}{\coqcode{PointSet_meet}{PointSet\_meet}}
\newcommand{\PointClassInteriorName}{\coq{PointClassInterior}}
\newcommand{\GBasicPointSetName}{\coqcode{GBasicPointSet}{GBasicPointSet}}
\newcommand{\GPointOpenName}{\coqcode{GPointOpen}{GPointOpen}}
\newcommand{\GPointSingleName}{\coqcode{GPointSingle}{GPointSingle}}
\newcommand{\GClosurePointName}{\coqcode{GClosurePoint}{GClosurePoint}}
\newcommand{\GPointNeighbourhoodName}{\coqcode{GPointNeighbourhood}{GPointNeighbourhood}}

\newcommand{\BallCenterBundleName}{\coqcode{BallCenterBundle}{BallCenterBundle}}

\newcommand{\CoqLibraryURL}{\url{https://github.com/mereoUnivers01/Tarski-Reconstruction}}

\definecolor{RSLRegionalInk}{RGB}{40,88,78}
\definecolor{RSLRegionalFill}{RGB}{232,244,239}
\definecolor{RSLPointInk}{RGB}{46,76,132}
\definecolor{RSLPointFill}{RGB}{232,238,250}
\definecolor{RSLBridgeInk}{RGB}{126,84,38}
\definecolor{RSLBridgeFill}{RGB}{251,239,220}
\definecolor{RSLArrow}{RGB}{73,82,96}

\tikzset{
  concept/.style={
    rectangle,
    rounded corners,
    draw,
    align=center,
    minimum height=1cm,
    minimum width=2.6cm,
    font=\small
  },
  arrow/.style={
    -{Latex[length=2mm]},
    thick
  },
  note/.style={
    font=\scriptsize,
    align=center
  }
}


\begin{document}

\title{From Regional Topology to Point-Class Topology in Tarski’s Geometry of Solids}
\author{Patrick Barlatier\\LISTIC, Universit{\'e} Savoie Mont Blanc\and Richard Dapoigny\\LISTIC, Universit{\'e} Savoie Mont Blanc}
\date{}
\maketitle

\begin{abstract}
Tarski's geometry of solids reconstructs point-like entities from concentric families of spherical regions rather than taking points as primitive. We formalize this reconstruction in Coq within a nominal mereological framework inspired by Le{\'s}niewski. Our main question is how the regional geometry of solids can induce a topology on reconstructed point-classes without adding new point individuals to the underlying mereological ontology.

We first develop the regional side of the construction. Three of Tarski's postulates concerning solids and their interior points, namely P2--P4, are derived as theorems of the mereological framework. Since Tarski's original interior-point definition allows concentric ball representatives to extend beyond the solid, we introduce the refined plural \texttt{sat\_InteriorPoint}. The mereological class generated by these representatives defines a regional interior operator \texttt{InteriorG}, for which we verify the four Kuratowski interior axioms. This regional construction also supports a boundary operation and an encoding of the RCC8 relations. It yields a topology of open solids, but it does not by itself provide a topology or a neighbourhood closure on reconstructed point-classes.

We therefore pass from regional solids to the setoid of ball representatives modulo \texttt{same\_center}. Each ball \(Q\) generates a stable basic point-plural \texttt{GBasicPointSet(Q)}, consisting of reconstructed point-classes having a representative regionally contained in \(Q\). These plurals form a basis for a metatheoretic topology on the setoid. The resulting topology is Hausdorff under Tarski's separation axiom \texttt{Three\_points} and is non-discrete under the explicit local richness hypothesis \texttt{BallCenterBundle}.

Finally, we define geometric neighbourhoods and the associated operator \texttt{GClosurePoint}. The central Coq theorem proves the four Kuratowski closure axioms; moreover, closure produces stable point-class plurals and respects extensional point-set equality. The lemma \texttt{GBasicPointSet\_as\_interior} identifies each ball-generated basic plural with the point-class plural induced by the regional part-plural of that ball. The formalization thus verifies the passage from a regional topology of solids to a Hausdorff topology and neighbourhood closure on reconstructed point-classes, without reifying those points as new mereological individuals.
\end{abstract}

\noindent\textbf{Keywords:} Tarski’s geometry of solids; Le{\'s}niewskian mereology; concentric point-classes; point-free topology; setoid quotients; Kuratowski closure; Coq formalization; nominal ontology

\medskip
\noindent\textbf{2025 Mathematics Subject Classification:} 03A05; 03B30; 03F25


\section{Introduction}\label{intro}
In the field of Qualitative Spatial Reasoning (QSR), mereotopology, the interplay between mereology and topology, is a common foundation for many theories \cite{Biacino1991,Roeper1997,Duntsch2004,Hahmann2013}. Theories of mereotopology follow the ideas of Whitehead’s point-free geometry \cite{WhiteHead1929} in which points are substituted by regions of space. Building on mereology either using a quasi-topological theory or by introducing topological notions to increase its expressiveness, e.g., in contact algebras or in the Region Connection Calculus (RCC8) \cite{Randell1992}, has gained significant attention in QSR research \cite{Menger1940,Varzi1994,Smith1996,Mormann2001}. However, these approaches often face challenges such as (i) the formal specification of topological boundaries (i.e., what kind of entity a boundary is, either a spatial object or an abstract entity?), (ii) the decidability of the resulting systems (the question of determining tractable fragments is still unanswered for many theories) \cite{Hahmann2012,Hahmann2013}, (iii) insufficient reasoning mechanisms without which the spatial representations are less useful and (iv) a shallow representation of mereology using a direct mapping of part of relations in first-order logic and resulting in poor expressiveness. To address these problems, we build on our previous work, i.e., the $\lambda$-MM open source library \cite{Barlat2023} which specifies an algebraic structure grounded in a syntactic, type-theoretic representation of mereology implemented in the Coq theorem prover.\footnote{\raggedright During the preparation of the formal development, Codex and ChatGPT Pro were used as exploratory aids for testing alternative formulations, exposing ambiguities, and suggesting possible proof patterns. All final definitions, theorem statements, Coq scripts, and mathematical interpretations were reviewed and accepted by the authors, who remain solely responsible for the paper.\par} The main objective is (i) to extend this library with an algebraic formulation of topological relations using appropriate primitives as an alternative to contact-based mereotopologies and (ii) to provide a unified theory incorporating mereology, geometry and topology. For this purpose, we will extend the $\lambda$-MM library implemented as a syntactic model grounded in \textit{Coq} \citep{Bertot2004,Coquand1988}. The library includes both mereology and Tarski’s geometry of solids within an inductive type space, all while maintaining the classical logic required by many mereological lemmas\footnote{\raggedright The Coq development is available at \CoqLibraryURL.\par}.

\section{Related Works}\label{related}
The present work belongs to the wide tradition of point-free and region-based geometry. In this tradition, regions or bodies are taken as primitive and points are reconstructed rather than assumed. Primarily, theories based on Whitehead's point-free geometry substitute points by regions of space as their primitive ontology \cite{WhiteHead1929}. Later formalizations have investigated the mereological and mereotopological structure of such systems. Tarski’s geometry of solids is especially relevant, and one of the closest developments of this seminal work is Gruszczyński and Pietruszczak’s full development of Tarski’s geometry of solids \cite{Gruszczynski2008}. Their assumed background is classical: solids are treated as mereological sums of balls, and point-like entities are recovered from region-based structure. By contrast, our construction separates two uses
of the data generated by boundary balls. In a first part, points which belong to the interior of a solid are saturated, while points in the boundary are not. The second alternative builds on concentric families of balls, which are
treated as equivalence classes of ball representatives under equality of concentric families. Both alternatives rely on $Gspace$, the space for solids and balls.
\par Another line of related work emerges from region-based topology and contact algebra (see e.g., \cite{Stell2000,Vakarelov2002,Duntsch2007}). Contact algebras treat regions as elements of a Boolean or distributive algebra equipped with a contact relation; local contact algebras and de Vries–Roeper–Dimov style dualities relate such structures to locally compact Hausdorff spaces. Dimov’s work on local contact algebras is a central contribution in this line \cite{Dimov2006}. Ivanova and Vakarelov's extended distributive contact lattices further weaken the Boolean setting by taking contact, dual contact, and nontangential inclusion as primitive notions. Smith and Varzi's distinction between bona fide and fiat boundaries also provides an important formal-ontological account of spatial demarcation \cite{SmithVarzi2000}; Smith's analysis of drawing lines on maps gives a complementary treatment of humanly induced spatial boundaries \cite{Smith1995}. The present work addresses a different but related problem: the construction of a verified neighbourhood and closure semantics for region-induced point-class plurals. These approaches operate in the extensional algebraic or formal-ontological register: regions are given as algebraic or spatial entities and contact-like or boundary relations are imposed between them. The proposed work operates in a preceding register. It reconstructs the ball representatives on which topology can be specified. The relevant data are not only regions and balls, but also equivalence classes of concentric ball representatives, stable point-class plurals, neighbourhoods, and the closure induced by them.
\par Recent work on contact join-semilattices and contact posets also weakens the algebraic basis of contact structures. Ivanova's contact join-semilattices \cite{Ivanova2019} drop Boolean complementation and reduce the algebraic structure needed for contact. Lipparini's contact posets \cite{Lipparini2023} retain only an order and a contact relation while obtaining strong model-theoretic properties. These developments show that region-based contact can be studied in weak algebraic settings. Our weakening is different: it is not a reduction of operations on already-given regions, but a reconstruction of the point-class plurals on which topological closure can be stated.
\par Pointless topology and formal topology also reject primitive points, but in a different way. In locale theory and pointless topology, the frame of opens is primitive and points are derived or optional. Constructive pointfree topology refines this idea through formal bases and cover relations; for example, Kawai and Sambin \cite{Kawai2019} define pointfree continuity by requiring operations on formal points to be induced by relations between bases. Grzegorczyk's point-free geometry is a further classical reconstruction in which point-like objects are represented by systems of regions \cite{Grzegorczyk1960}. Our construction is point-free, but in a Tarskian sense: we do not start from a frame of opens, but from balls, solids, and the quotient structure induced by concentricity. 
\par Existing point-free, contact-algebraic, and nominal techniques each capture part of the landscape. What is new here is their combination in a Tarskian mereotopology where points are reconstructed from balls and then organized into a genuine point-open topology. The contribution is not the introduction of setoids or quotients as such. Rather, the paper identifies the quotient structure needed to make Tarski's reconstruction of points support a Kuratowski closure operator in a nominal mereological formalization. The formal mechanism is not a simple addition of a boundary region to a body. It proceeds by passing from balls to ball representatives, identifying them modulo the relation of having the same center, defining point-class plurals of representatives, and finally defining closure
by neighbourhood intersection.
\par The second part explains how concentric point-classes are closed under topological adherence. The gap addressed in the present paper is therefore not the absence of point-free accounts of space or accounts based on regions. Rather, it is the lack of a formal bridge between Tarski’s reconstruction of points from concentric balls and a Kuratowski closure operator acting on the corresponding point-classes. The present work supplies this bridge by treating ball representatives through setoids and by defining closure through point-open neighbourhoods.
\section{The nominal mereological background}
The $\lambda$-MM library \cite{Barlat2023} relies on Le{\'s}niewski's stratified systems \cite{Lesniewski1916}. We recall only the nominal and mereological notation used later for the point-class construction.
\subsection{Notation and mathematical dependencies}

For readability, Table~\ref{tab:notation-guide} collects the notations
that carry the main construction. The table is not intended to replace the
formal definitions given below; it fixes the informal reading of the
symbols before the technical development begins. Throughout the paper,
uppercase letters \(P,Q,R,\ldots\) denote individual names or ball
representatives, whereas lowercase letters \(a,b,u,v,\ldots\) denote
plural names. Thus \(\eta\,P\,a\) is read as singular inclusion of the
individual \(P\) under the plural name \(a\).

\begin{center}
\begin{minipage}{\textwidth}
\refstepcounter{table}\label{tab:notation-guide}
\textbf{Table \thetable: Notational guide for the point-class construction}\\[2pt]
\footnotesize
\setlength{\tabcolsep}{4pt}
\begin{tabular}{@{}p{0.25\textwidth}p{0.25\textwidth}p{0.42\textwidth}@{}}
\hline
Symbol & Meaning & Informal interpretation \\
\hline
\(\eta\,P\,a\) & nominal singular inclusion & the individual name or representative \(P\) falls under the plural name \(a\) \\
\(\coq{klass}(a)\) & generated m-class & the mereological individual generated by the plural \(a\) \\
\(\coq{pt}(Q)\) & part-plural & the plural of parts of \(Q\) \\
\(\coq{Point}(Q)\) & point-class generated by \(Q\) & the plural of balls concentric with the representative \(Q\) \\
\(\SameCenter(P,Q)\) & representative equivalence & \(P\) and \(Q\) determine the same concentric point-class \\
\(\coq{InteriorG}\) & regional interior & the interior operation in the regional topology of solids \\
\(\coq{InteriorPoint}\) & Tarskian interior-point predicate & the regional predicate used in Tarski's reconstruction of points \\
\(\PointClassInteriorName\) & induced point-class interior & the point-class plural induced by regional data \\
\(\GBasicPointSetName\) & basic point-plural & point-class plural generated by a regional ball \\
\(\GPointOpenName\) & geometric point-open plural & open plural generated by the ball basis \\
\(\GClosurePointName\) & Kuratowski closure & neighbourhood closure on stable point-class plurals \\
\hline
\end{tabular}
\end{minipage}
\end{center}

\paragraph{Mathematical dependencies.}
This paper uses the $\lambda$-MM library through a small number of
ingredients: nominal plurals and the membership relation \(\eta\);
mereological classes and parthood, expressed by \(\coq{klass}\) and
\(\coq{pt}\); the reconstruction of point-like objects from concentric
balls; representative equivalence under \SameCenter{}; the regional
interior structure inherited from Tarski's geometry of solids; and the
induced point-class operations \PointClassInteriorName{} and
\GBasicPointSetName{} and \GClosurePointName{}. The explicitly geometric assumptions used by the
library are recalled later in Subsection~\ref{subsec:formal-geometric-background}.
Thus the reader only needs to keep the following dependency chain in view:

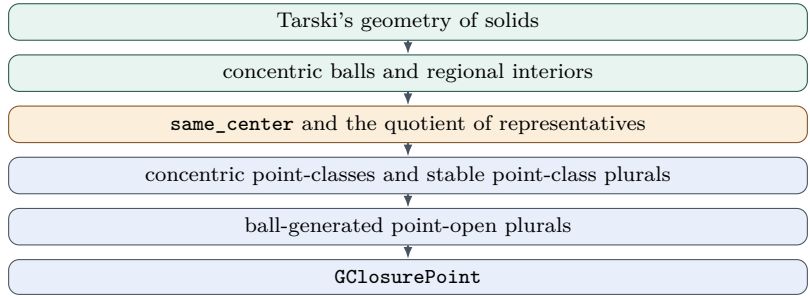
\begin{center}
\small
\begin{tikzpicture}[
  depnode/.style={
    rectangle,
    rounded corners,
    draw=RSLArrow,
    fill=RSLPointFill,
    align=center,
    minimum width=0.64\textwidth,
    minimum height=0.46cm,
    font=\scriptsize
  },
  deparrow/.style={-{Latex[length=1.6mm]}, thick, draw=RSLArrow},
  node distance=0.18cm
]
\node[depnode, fill=RSLRegionalFill, draw=RSLRegionalInk] (g) {Tarski's geometry of solids};
\node[depnode, below=of g, fill=RSLRegionalFill, draw=RSLRegionalInk] (b) {concentric balls and regional interiors};
\node[depnode, below=of b, fill=RSLBridgeFill, draw=RSLBridgeInk] (s) {\SameCenter{} and the quotient of representatives};
\node[depnode, below=of s] (p) {concentric point-classes and stable point-class plurals};
\node[depnode, below=of p] (o) {ball-generated point-open plurals};
\node[depnode, below=of o] (c) {\GClosurePointName{}};
\draw[deparrow] (g) -- (b);
\draw[deparrow] (b) -- (s);
\draw[deparrow] (s) -- (p);
\draw[deparrow] (p) -- (o);
\draw[deparrow] (o) -- (c);
\end{tikzpicture}
\end{center}

\subsection{From Le{\'s}niewski's Ontology to $\lambda$-MM}
Le{\'s}niewski's Ontology (LO), or logic of names, manipulates names rather than predicates as its basic syntactic material. Its central relation is the binary functor $\eta$: the sentence $\eta \: A \: x$ says that the name \(A\) falls under the classification designated by the name \(x\). This is not set-theoretic membership, but an internal nominal relation. For example, if \(A\) designates Nicole and \(B\) designates ``the neighbour's daughter'', then \(\eta B A\) expresses the corresponding nominal identification inside the language, rather than by appealing to an external assignment. The formal role of $\eta$ is close to the ``:'' notation of type theory \cite{Simons1987}, but without separating objects, types, and interpretation into distinct meta-theoretic layers.
\par In $\lambda$-MM, Coq's dependent type theory provides a syntactic model of this nominal framework. Names are objects of an inductive type \(N\), and $\eta$ is represented as a judgmental relation of type \(N \rightarrow N \rightarrow Prop\). We distinguish simple relations, of type \(N \rightarrow N \rightarrow Prop\), from compound relations obtained by composing $\eta$ with unary or binary name constructors. For a constructor \(\phi\), membership of an object \(X\) in the name \((\phi a)\) is inferentially equivalent to singular inclusion of its singleton name \(\iota X\):
\begin{lemma}\label{In_to_eta} 
$\forall \: a : N , \: \forall \: X : object, In \; (\phi \: a) \: X \: \leftrightarrow \: \eta \: (\iota \: X) \: (\phi \: a)$	
\end{lemma}
in which $\iota$ denotes the singleton constructor\footnote{It maps every object $X$ to its associated name, $\iota \: X$.}. The basic relations used below are listed in Table~\ref{table1}.

\begin{table}[t]
\caption{Basic relations between names in $\lambda$-MM\label{table1}}

\begin{tabular}[!h]{ c }
\textit{Simple relations} \\
\begin{tabular}{ c  c  c  c  c }
\hline
$\eta$ & singular inclusion & $\forall P \; Q, \eta \: P \: Q \;$ & &  \\
$\equiv$ & equality of individual names & $\forall P \; Q, P \: \equiv \: Q $ &$\leftrightarrow$& $\; \eta \: P \; Q \wedge \eta \; Q \; P$ \\
$\approx$ & equality of plural names & $\forall a \: b, a \approx b$ &$\leftrightarrow$& $\forall P, \eta \: P \: a \leftrightarrow \; \eta P \: b$ \\
$\subseteq$ & name inclusion & $\forall a \: b, a \subseteq b$ &$\leftrightarrow$& $\; \forall P, \; \eta \: P \: a \rightarrow \; \eta P \: b$ \\
\hline
\end{tabular} \\
\textit{Compound relations }\\
\begin{tabular}{ c  c  c  c  c }
\hline
$\mathit{neg}$ & negation of name & $\forall P \; Q, \eta \: P \; (neg \; Q)$ &$\leftrightarrow$& $(\eta \; P \; P \wedge \neg(\eta \; P \; Q))$ \\
$\cap$ & name conjunction & $\forall a \: b, a \cap b$ &$\leftrightarrow$& $ \forall P, \; (\eta \: P \: a \wedge \eta \: P \: b)$ \\
$\cup$ & name disjunction & $\forall a \: b, a \cup b$ &$\leftrightarrow$& $\; \forall P, \; (\eta \: P \: a \vee \eta \: P \: b) $ \\
\hline
\end{tabular}
\end{tabular}
\end{table}
Two constant plural names are specified in LO, which are respectively the empty name ($\Lambda$) defined as the contradictory name, and the universal name ($V$). 
\subsection{From Le{\'s}niewski's Mereology to $\lambda$-MM}
Le{\'s}niewski's mereology treats a class not as a set of members but as a concrete whole, or mereological sum, composed of parts. Thus an m--class is itself an individual, and its elements are understood through the part--whole relation. In $\lambda$-MM, \(\eta A (klass\,a)\) says that \(A\) is the individual generated by the plural name \(a\), whereas \(\eta A a\) treats \(A\) as falling directly under the plural \(a\). M--classes are extensional: they are determined by their parts \cite{Asenjo1977}. The relevant constructors, such as \(pt\), \(ppt\), and \(klass\), are all compound nominal relations satisfying Lemma~\ref{In_to_eta}. The definition of \(klass\) is recalled because it is later used in the regional reconstruction theorem:
\begin{definition}[\coq{m-class}]\label{mklass}\mbox{}\\
$\forall A \: a, \: \eta \: A \: (klass \: a) \triangleq (Individual \: A \: \wedge \; (\forall \: B, \eta \: B \: a \rightarrow \eta \: B \: (pt \: A)) \: \wedge \;
(\forall \: B, \eta \: B (pt \: A) \rightarrow $\\ $\exists \: C \: D, \eta \: C \: a \wedge \eta \: D (pt \: C) \wedge \eta \: D (pt \: B)))$.
\end{definition}
The term \(Individual\,A\) is equivalent to \(\eta A A\). The full development in \cite{Barlat2023} yields a quasi-Boolean algebra\footnote{A Boolean algebra without a bottom element.} when the order is identified with the \(pt\) relation \cite{Tarski1956a,Clay1974}. Thus \(A\leq B\) and \(\eta A (pt\,B)\) are two presentations of the same parthood structure: the former order-theoretic, the latter internal to the nominal language.
\section{The regular-open background of $\lambda$-MM}\label{TM}
\subsection{Introduction of Meet and Join Operators}
In $\lambda$-MM, the whole space of mereology is defined as an individual name called $Universe$ such as:
\begin{definition}[\coq{Universe}]\label{Univ}\mbox{}\\
$\forall \: P \: : N, \; \eta \: P \: Universe \triangleq \eta \: P \: (klass \: V).$
\end{definition}
The Boolean presentation of $\lambda$-MM uses the part-of relation \(pt\) as its order. We only recall the induced join and meet operations, since they provide the regular-open background used later:
\begin{definition}[\coq{join}]\label{join}\mbox{}\\
$\forall P \: Q \: R : N, \; \eta \: R \: (b\_sum \: P \: Q) \; \triangleq \; (\eta \: P \: P \; \wedge \; \eta \: Q \: Q \; \wedge \; \eta \: R \: (klass \: ((pt \: P) \cup (pt \: Q)))).$
\end{definition}
Since no bottom element exists, meet is used locally when the conjunction of parts is defined:
\begin{definition}[\coq{meet}]\label{meet}\mbox{}\\
$\forall P \: Q \: R : N, \; \eta \: R \: (b\_prod \: P \: Q) \; \triangleq \; (\eta \: P \: P \; \wedge \; \eta \: Q \: Q \; \wedge \; \eta \: R \: (klass \: ((pt \: P) \cap (pt \: Q)))).$
\end{definition}
These operators satisfy the expected uniqueness and complement properties in the library.
\subsection{Mereology as a topology of individual names}\label{Btopo}
The regular-open interpretation developed in \cite{Dapoigny2026} views mereological individuals as non-empty open subsets of a topology \(\langle N,\mathbb{T}_{MM}\rangle\). In that setting, the interior of an individual \(Q\) is \(Q\) itself:
 \begin{definition}[\coq{interior}]\label{interior}\mbox{}\\
$\forall \: P \: Q \: : N, \; \eta \: P \: (interior \: Q) \triangleq \: \eta \: Q \: Q \: \wedge \: \eta \: P \: Q$
\end{definition}
With this definition, inclusion, intersection, union, and complement correspond respectively to \(pt\), \(b\_product\), \(b\_sum\), and \(compl\), following the regular-open reading of mereology \cite{Mormann2001}. The relevant axioms are the interior-dual form of the Kuratowski axioms \cite{Pervin1964}:
\begin{enumerate}
\item Preservation of the total space: $int(\mathbb{T}_{MM}) = \mathbb{T}_{MM}$.
\item Intensiveness: $\forall A \subseteq \mathbb{T}_{MM}: int(A) \subseteq A$.
\item Idempotence: $\forall A \subseteq \mathbb{T}_{MM}: int(int A)) = (int  A)$.
\item Preservation of binary intersections: $\forall A, B \subseteq \mathbb{T}_{MM}: int(A \cap B) = int(A) \cap int(B)$
\end{enumerate}
The corresponding Coq theorems are recalled only to fix notation:
\begin{theorem}[\coq{int\_of\_space}]\label{intsp}\mbox{}
\begin{coqdoccode}
(interior Universe) $\approx$ Universe.
\end{coqdoccode}	
\end{theorem}
\begin{theorem}[\coq{int\_intensive}]\label{intint}\mbox{}
\begin{coqdoccode}
\ensuremath{\forall} Q, (interior Q) $\subseteq$ Q.
\end{coqdoccode}	
\end{theorem}
\begin{theorem}[\coq{int\_open\_intensive}]\label{intopint}\mbox{}
\begin{coqdoccode}
\ensuremath{\forall} Q, $\eta$ Q Q $\rightarrow$ (interior Q) $\approx$ Q.
\end{coqdoccode}	
\end{theorem}
\begin{theorem}[\coq{int\_idempotent}]\label{intidem}\mbox{}
\begin{coqdoccode}
\ensuremath{\forall} Q, (interior (interior Q)) $\approx$ (interior Q)
\end{coqdoccode}	
\end{theorem}
\begin{theorem}[\coq{bin\_intersections}]\label{intinter}\mbox{}
\begin{coqdoccode}
\ensuremath{\forall} \coqdocvar{P} \coqdocvar{Q} \coqdocvar{R} \coqdocvar{S}, \coqdocdefinition{$\eta$} \coqdocvar{P} \coqdocvar{P} \ensuremath{\land} \coqdocdefinition{$\eta$} \coqdocvar{Q} \coqdocvar{Q} \ensuremath{\land} \coqdocdefinition{$\eta$} \coqdocvar{R} (interior \coqdocvar{P}) \ensuremath{\land} \coqdoceol
\coqdocindent{2.00em}
\coqdocdefinition{$\eta$} S (interior \coqdocvar{Q}) $\rightarrow$ (\coqdocvar{b\_product} \coqdocvar{R} \coqdocvar{S}) $\approx$ interior (\coqdocvar{b\_product} \coqdocvar{P} \coqdocvar{Q}).
\end{coqdoccode}	
\end{theorem}
Closure is then introduced as the complement of the interior of the complement:
\begin{definition}[\coq{closure}]\label{closure}\mbox{}\\
$\forall \: P \: Q, \: \eta \: P \: (closure \: Q) \triangleq \; (\eta \: Q \: Q \wedge \; (exists \: R \: S, \; \eta \: R \: (compl \: Q) \: \wedge \: \eta \: S \: (interior \: R) \: \wedge $ \\ $ \eta \: P \: (compl \: S))).$	
\end{definition}
Regularity states that the interior of the closure of \(Q\) coincides with \(Q\):
\begin{lemma}[\coq{regular}]\label{regul}\mbox{}
\begin{coqdoccode}
\ensuremath{\forall} \coqdocvar{P} \coqdocvar{Q}, (\coqdocdefinition{$\eta$} \coqdocvar{Q} \coqdocvar{Q} \ensuremath{\land} \coqdocdefinition{$\eta$} \coqdocvar{P} (\coqdocvar{closure} \coqdocvar{Q}) \ensuremath{\rightarrow} \coqdocvar{Q} $\equiv$ (\coqdocvar{interior} \coqdocvar{P})).\coqdoceol
\end{coqdoccode}	
\end{lemma}
The mereological boundary of \(Q\) is defined as the intersection of the closure of \(Q\) and the closure of its complement:
\begin{definition}[\coq{boundary}]\label{bound0}\mbox{}\\
$\forall \: P \: Q \: : N, \: \eta \: P \: (boundary\: Q) \: \triangleq \: (\eta \: Q \: Q \; \wedge \; \eta \: P \: ((closure \: Q) \: \cap \: (closure \: (compl \: Q)))).$
\end{definition}			
By definition \ref{bound0}, we prove that this boundary is empty and then, that the topology is clopen:
\begin{theorem}[\coq{clopen}]\label{clopen}\mbox{}
\begin{coqdoccode}
\ensuremath{\forall} \coqdocvar{Q}, $\eta$ \coqdocvar{Q} \coqdocvar{Q} \ensuremath{\land} \ensuremath{\lnot} (\coqdocvar{Q} $\equiv$ \coqdocvar{Universe}) \ensuremath{\rightarrow} \coqdocvar{boundary} \coqdocvar{Q} $\approx$ $\Lambda$.\coqdoceol
\end{coqdoccode}	
\end{theorem}
This individual-name topology is useful as background but too coarse for the point-class closure developed below. The next section recalls only the Tarskian structures needed to obtain the regional source of point-class data.
\section{Topological Interpretation of Tarski's Geometry of Solids}\label{TopoTarski}
In the second part of \cite{Dapoigny2026}, $\lambda$-MM was extended with Tarski's geometry of solids, using balls as primitive regional objects. Tarski's system is naturally connected with regular open sets and with Euclidean models of \(\mathbb{R}^n\) \cite{Tarski1956b,Bennett2001}. Here we use only the part of that development needed for the present paper: balls, solids, point-like families of concentric balls, and the regional topology from which point-class data will be induced. The key limitation of the previous regional presentation is not its regional topology itself, but the fact that Kuratowski closure on reconstructed locations requires a point-class carrier rather than a mereological solid.
\subsection{Basic Tarskian structures}
\label{sec:tarskian-structures}
Tarski's geometry has been extended on the basis of mereology in a previous work \cite{Dapoigny2026}. The previous subsection recalled that the nominal mereological framework of \(\lambda\)-MM admits a regular-open interpretation. It is now refined in the specific setting of Tarski's geometry of solids. On the one hand, the model of Mereology coincides with that of a boolean algebra as a lattice structure without zero (see e.g., \cite{Tarski1956a,Clay1974,Barlat2023}). On the other hand, Tarski has given a model for atomless Boolean algebra, consisting of the family of regular open sets of an Euclidean space and the relation of set-inclusion \cite{Tarski1956b}. Besides balls, Tarski's geometry includes solids as sets of balls and points as infinite sets of concentric balls. The key point is that the whole space should not be treated as an homogeneous region and we need to distinguish solids as sets of balls from sets of points as infinite structures generated by a given ball. For that purpose, we introduce a regional structure that includes solids and balls as specific solids \(\coq{Gspace}\).
\par In the Coq development, the name of all balls is denoted by \(\coq{balls}\), and the geometric space is obtained as the mereological class of balls:
\begin{definition}[\coq{Gspace}]\label{gspace}\mbox{}\\
$\forall \: P \: : N, \; \eta \: P \: \coq{Gspace} \triangleq \: \eta \: P \: (klass \: balls)$
\end{definition}
Note that since \(\coq{Gspace}\) is specified as a mereological class, it is an individual name and therefore satisfies the uniqueness property. A solid is then understood as an individual generated from a suitable plurality of balls. Thus, the regional structure of the theory is body-based: balls and solids are the primitive geometric material. It results that any ball is a member of \(\coq{Gspace}\):
\begin{lemma}[\coq{in\_Gspace}]\label{inGspace}\mbox{}
\begin{coqdoccode}
\ensuremath{\forall} \coqdocvar{A}, \coqdocdefinition{$\eta$} \coqdocvar{A} \coqdocvar{balls} \ensuremath{\rightarrow} \coqdocdefinition{$\eta$} \coqdocvar{A} (\coqdocdefinition{pt} \: \coqdocvar{Gspace}).
\end{coqdoccode}
\end{lemma}
\par which is interpreted as $A \: \in \: balls \rightarrow A \: \in \: pt(Gspace)$.
\par The point-set structure first requires the definition of \coq{Point}. Skipping details, points rely on the concentricity introduced by Tarski as an equivalence relation. They are specified as the set of all balls which are concentric with a given one.
\begin{definition}[\coq{Point}]\label{point}\mbox{}
$\forall \: P \: Q \: : N, \; \eta \: P \: (\coq{Point}\,Q) \triangleq \: (\eta \: P \: balls \: \wedge \:\eta \: Q \: balls \: \wedge \: \eta \: P \: (\coq{Concent}\,Q))$
\end{definition}
Notice that \((\coq{Concent}\,Q)\) is a plural, that is, a set of all balls concentric with ball \(Q\) and \(P\) is one of these balls. While balls denote first-order objects, points as set of balls are second-order objects and fall within the scope of a monadic second-order theory. Since concentricity is an equivalence relation, \coq{Point} inherits reflexivity, symmetry and transitivity. If \(P\) is a ball, then a point as the set of concentric balls with \(P\) exists:
\begin{lemma}[\coq{ex\_point}]\label{expt}\mbox{}
\begin{coqdoccode}
\ensuremath{\forall} \coqdocvar{P}, \coqdocdefinition{$\eta$} \coqdocvar{P} \coqdocvar{balls} \ensuremath{\rightarrow} \ensuremath{\exists} \coqdocvar{Y}, \coqdocdefinition{$\eta$} \coqdocvar{Y} (\coq{Point} \coqdocvar{P}).
\end{coqdoccode}
\end{lemma}
From the point definition, two further facts are worth detailing: points as
parts and equality of point plurals induced by concentric representatives.
\begin{lemma}[\coq{point\_as\_parts}]\label{pointsParts}\mbox{}
\begin{coqdoccode}
\ensuremath{\forall} \coqdocvar{P} \coqdocvar{B}, \coqdocdefinition{$\eta$} \coqdocvar{B} (\coq{Point} \coqdocvar{P}) \ensuremath{\rightarrow} \coqdocdefinition{$\eta$} \coqdocvar{B} (\coqdocdefinition{pt} \coqdocvar{P}) $\vee$ \coqdocdefinition{$\eta$} \coqdocvar{P} (\coqdocdefinition{pt} \coqdocvar{B}).
\end{coqdoccode}
\end{lemma}
\begin{theorem}[\coq{point\_sets}]\label{pointsSets}\mbox{}
\begin{coqdoccode}
\ensuremath{\forall} \coqdocvar{P} \coqdocvar{Q}, \coqdocdefinition{$\eta$} \coqdocvar{P} (\coq{Point} \coqdocvar{Q}) \ensuremath{\rightarrow} (\coq{Point} \coqdocvar{P}) $\approx$ (\coq{Point} \coqdocvar{Q}).
\end{coqdoccode}
\end{theorem}
Thus any ball representative \(P\) belonging to the point generated by
\(Q\) determines the same point-class as \(Q\). This equivalence theorem is
the formal reason why the later quotient relation can identify ball
representatives through equality of their \coq{Point} plurals.
Another important geometric relation can be introduced: betweenness between
points. Given two points, a third point lies between them iff there exist
three balls concentric with them which are in the external diametric relation:
\begin{definition}[\coq{Betweenness}]\label{betw}\mbox{}\\
$\forall A \: B \: C, \; \eta \: A \: (\coq{btw}\,B \: C) \triangleq
(\eta \: A \: balls \wedge \eta \: B \: balls \wedge \eta \: C \: balls \wedge$\\
$\exists \: D \: E \: F, \: \eta \: D \: (\coq{Point}\,B) \wedge
\eta \: E \: (\coq{Point}\,A) \wedge \eta \: F \: (\coq{Point}\,C) \wedge
\eta \: D \: (\coq{ED}\,E \: F)).$
\end{definition}
Co-concentric balls \(P\) and \(Q\) give rise to a common set. The next definition following Tarski's idea refers to the solid considered as a set of balls:
\begin{definition}[\coq{Solid}]\label{solid}\mbox{}\\
$\forall \: P \: : N, \: \eta \: P \: Solid \triangleq (Individual \: P \wedge \exists \: b, \: b \subseteq balls \: \wedge \; \eta \: P \: (klass \: b)).$
\end{definition}
A similar definition has been given in \cite{Clay2021}. It is more precise than Tarski's description since it involves a particular subset of balls giving rise to a corresponding m--class. However, Tarski has pointed out that the concept of ``sum'' should be understood as collection in Le{\'s}niewski's formalism. Precisely, using definition \ref{solid} we can derive the following equivalence which fully agrees with Tarski's view:
\begin{lemma}[\coq{solid\_as\_coll}]\label{solColl}\mbox{}
\begin{coqdoccode}
\coqdockw{\ensuremath{\forall}} \coqdocvar{P}, \coqdocdefinition{$\eta$} \coqdocvar{P} \coqdocvar{Region} \ensuremath{\leftrightarrow} \coqdocdefinition{$\eta$} \coqdocvar{P} (\coqdocdefinition{coll} \coqdocvar{balls}).
\end{coqdoccode}
\end{lemma}
Since a region is a set of balls and provided that balls are open sets, their union is open as well. Lemmas \ref{ballsol} and \ref{spacesol} show respectively that any ball is a solid (the converse is not true), and that the geometric space is also a solid.
\begin{lemma}[\coq{any\_ball\_is\_a\_solid}]\label{ballsol}\mbox{}
\begin{coqdoccode}
\coqdockw{\ensuremath{\forall}} \coqdocvar{P}, \coqdocdefinition{$\eta$} \coqdocvar{P} \coqdocvar{balls} \ensuremath{\rightarrow} \coqdocdefinition{$\eta$} \coqdocvar{P} \coqdocvar{Solid}.
\end{coqdoccode}
\end{lemma}
\begin{lemma}[\coq{whole\_space\_is\_a\_solid}]\label{spacesol}\mbox{}
\begin{coqdoccode}
\coqdockw{\ensuremath{\forall}} \coqdocvar{P}, \coqdocdefinition{$\eta$} \coqdocvar{P} \coqdocvar{Gspace} \ensuremath{\rightarrow} \coqdocdefinition{$\eta$} \coqdocvar{P} \coqdocvar{Solid}.
\end{coqdoccode}
\end{lemma}
\par Skipping details, we now focus on the last definition proposed by Tarski. The point \(a\) is an interior point of the solid \(Q\) if there exists a ball \(P'\) which is at the same time an element of the point \(a\) and a part of the region \(Q\). It has the following form:
\begin{definition}[\coq{InteriorPoint}]\label{intpt}\mbox{}\\
$\forall \: P \: Q : N, \: \eta \: P \: (\coq{InteriorPoint}\,Q) \triangleq (\eta \: Q \: Solid \wedge \exists \: P', \: \eta \: P \: (\coq{Point}\,P') \: \wedge \; \eta \: P' \: (pt \: Q)).$
\end{definition}
When Tarski refers to the point \(a\), he involves a set of all concentric balls to a given ball; therefore \coq{Point}\,\(P'\) plays here the role of \(a\) and \(P'\) is one of its balls. It follows that \(P\) is also one of these balls, provided that ball \(P'\) is a part of region \(Q\). Thus, points are seen as filters of regions. In summary the definition of the interior point describes a set which in Le{\'s}niewski's formalism can be related to a m--class. It is important to see that such a definition asserts that while all the referred points belong to the solid, their related concentric balls expand outside the region. This aspect is crucial for the definition of the topological interior of a region.
\par Beside lemmas proving that for any solid, there exists at least an interior point and that any individual point is a ball we can formalize Tarski's postulates 2 and 3 as theorems. Since an interior point corresponds to a set, the sum of this set is an individual according to mereological principles and then, it refers to a regular open set:
\begin{theorem}[\coq{Tarski\_P2}]\label{TP2}\mbox{}
\begin{coqdoccode}
\coqdockw{\ensuremath{\forall}} \coqdocvar{Q}, \coqdocdefinition{$\eta$} \coqdocvar{Q} \coqdocdefinition{Solid} \ensuremath{\rightarrow} \ensuremath{\exists} \coqdocvar{R}, \coqdocdefinition{$\eta$} \coqdocvar{R} (\coqdocdefinition{klass} (\coqdocdefinition{InteriorPoint} \coqdocvar{Q})).
\end{coqdoccode}
\end{theorem}
\begin{theorem}[\coq{Tarski\_P3}]\label{TP3}\mbox{}
\begin{coqdoccode}
\coqdockw{\ensuremath{\forall}} \coqdocvar{R} \coqdocvar{Q}, \coqdocdefinition{$\eta$} \coqdocvar{R} (\coqdocdefinition{klass} (\coqdocdefinition{InteriorPoint} \coqdocvar{Q})) \ensuremath{\rightarrow} \coqdocdefinition{$\eta$} \coqdocvar{R} \coqdocdefinition{Solid}.
\end{coqdoccode}
\end{theorem}
Theorem \ref*{TP2} states that if \(Q\) is a region, the m--class of all its interior points is non-empty regular open set. This results from the properties of any m--class. Alternatively, in theorem \ref*{TP3} if we have a class of points that is a regular open set (i.e., a m--class), then it corresponds to an individual \(R\) which is a solid.
\par The interpretation of postulate 4 demonstrates that the inclusion of sets of interior points for two solids implies that their respective open sets are in the part-of relation:
\begin{theorem}[\coq{Tarski\_P4}]\label{TP4}\mbox{}
\begin{coqdoccode}
\coqdoceol
\coqdocindent{2.00em}\ensuremath{\forall} \coqdocvar{P} \coqdocvar{Q} \coqdocvar{R} \coqdocvar{S}, \coqdocdefinition{$\eta$} \coqdocvar{P} \coqdocdefinition{Solid} \ensuremath{\land} \coqdocdefinition{$\eta$} \coqdocvar{Q} \coqdocdefinition{Solid} \ensuremath{\land} (\coqdocdefinition{InteriorPoint} \coqdocvar{P}) $\subseteq$ (\coqdocdefinition{InteriorPoint} \coqdocvar{Q}) \ensuremath{\land} \coqdoceol
\coqdocindent{2.00em}\coqdocdefinition{$\eta$} \coqdocvar{R} (\coqdocdefinition{klass} (\coqdocdefinition{InteriorPoint} \coqdocvar{P})) \ensuremath{\land}
\coqdocdefinition{$\eta$} \coqdocvar{S} (\coqdocdefinition{klass} (\coqdocdefinition{InteriorPoint} \coqdocvar{Q})) \ensuremath{\rightarrow} \coqdocdefinition{$\eta$} \coqdocvar{R} (\coqdocdefinition{pt} \coqdocvar{S}).
\end{coqdoccode}
\end{theorem}
\subsection{Interior and saturated interior points}
As underlined in Section~\ref{TopoTarski}, the original specification of interior points is not sufficient to construct the base for a topology. Therefore, we have slightly modified Tarski's definition \ref{intpt} to constrain balls of interior points to be saturated, i.e., constrained to be inside the solid:
\begin{definition}[\coq{SatInteriorPoint}]\label{satintpt}\mbox{}\\
$\forall \: P \: Q, \: \eta \: P \: (\coq{sat\_InteriorPoint}\,Q) \triangleq (\eta \: Q \: Solid \: \wedge \; \exists \: P', \: \eta \: P \: (\coq{Point}\,P') \: \wedge \: \eta \: P' \: (pt \: Q) \: \wedge \: (\eta \: P' \: (pt \: P) \: \rightarrow$\\ $(P \: \equiv \: P'))).$
\end{definition}
As expected, in Tarski's definition of \coq{InteriorPoint}, ball \(P\) can extend outside of solid \(Q\), since it is only its center which is constrained to be inside \(Q\). By contrast, definition \ref{satintpt} constrains ball \(P\) to be inside \(Q\). Using lemma \ref{pointsParts}, the term $\eta P \; (\coq{Point}\,Y)$ yields two cases: $\eta \: Y \; (pt \: P)$ or $\eta \: P \; (pt \: Y)$. Then we consider two terms, $\eta \: Y \; (pt \: Q)$ and the condition: $(\eta \: Y \: (pt \: P) \: \rightarrow \: (P \equiv Y)$. Let us consider the first alternative, i.e., $\eta \: Y \: (pt \: P)$. Applying the condition to this term provides $(P \equiv Y)$ and using the remaining term $\eta \: Y \: (pt \: Q)$ it follows that $\eta \: P (pt \: Q)$. The second alternative, $\eta \: P \: (pt \: Y)$ can be merged with the term $\eta \: Y \: (pt \: Q)$ by applying transitivity of the part-of relation and then deriving easily $\eta \: P \: (pt \: Q)$. It follows that any ball \(Y\) concentric with \(P\) guarantees that \(P\) always remains inside \(Q\).
\par All lemmas derived with Tarski's definition \ref{intpt} remain valid but we can prove specific lemmas. For example, (i) the geometric space is an open set built with the m--class of its interior points and (ii) if \(Q\) is an open solid, then it coincides with the m--class of its interior points. The term \coq{SatInteriorPoint} refers to a plural name whose members are balls such as \(P\). In mereology, each plural \(a\) can be the generator of a class \(A\) such that $\eta \: A \: (klass \: a)$. Alternatively, a mereological lemma states that: $a \: \subseteq \: (pt \: A)$. It means that each element of the plural is also a part of the related class. Therefore, a mereological class can be seen as a generalization of its related plural. We introduce the topological interior as the generalization of the set of saturated interior points:
\begin{definition}[\coq{InteriorG}]\label{intgeom}\mbox{}\\
$\forall \: P \: Q \: : N, \: \eta \: P \: (\coq{InteriorG}\,Q) \triangleq \; \eta \: P \: (klass \: (\coq{sat\_InteriorPoint}\,Q)).$
\end{definition}
that is, \(P\) is the topological interior of any region \(Q\) iff it is the m--class of its saturated interior points. It is a singular name that refers to all saturated interior points. All Kuratowski's axioms of subsection \ref{Btopo} are rewritten as follows:
\begin{theorem}[\coq{int\_of\_space}]\label{intGsp}\mbox{}
\begin{coqdoccode}
(InteriorG Gspace) $\equiv$ Gspace.
\end{coqdoccode}
\end{theorem}
\begin{theorem}[\coq{int\_intensive}]\label{intGint}\mbox{}
\begin{coqdoccode}
\ensuremath{\forall} Q, (InteriorG Q) $\subseteq$ Q.
\end{coqdoccode}
\end{theorem}
\begin{theorem}[\coq{int\_open\_intensive}]\label{intopGint}\mbox{}
\begin{coqdoccode}
\ensuremath{\forall} Q, $\eta$ Q Solid $\rightarrow$ (InteriorG Q) $\equiv$ Q.
\end{coqdoccode}
\end{theorem}
\begin{theorem}[\coq{int\_idempotent}]\label{intGidem}\mbox{}
\begin{coqdoccode}
\ensuremath{\forall} Q, $\eta$ Q Solid $\rightarrow$ \coqdoceol
\coqdocindent{2.00em} (InteriorG (InteriorG Q)) $\equiv$ (interiorG Q)
\end{coqdoccode}
\end{theorem}
\begin{theorem}[\coq{bin\_intersections}]\label{intGinter}\mbox{}
\begin{coqdoccode}
\ensuremath{\forall} \coqdocvar{P} \coqdocvar{Q}, \coqdocdefinition{$\eta$} \coqdocvar{P} \coqdocvar{Solid} \ensuremath{\land} \coqdocdefinition{$\eta$} \coqdocvar{Q} \coqdocvar{Solid} $\rightarrow$ \ensuremath{\forall} \coqdocvar{R} \coqdocvar{S} \coqdocvar{Z}, \coqdoceol
\coqdocindent{2.00em} \coqdocvar{R} $\equiv$ (InteriorG P) \ensuremath{\land} 
\coqdocvar{S} $\equiv$ (InteriorG Q) \ensuremath{\land} \coqdocvar{Z} $\equiv$ (\coqdocvar{b\_product} \coqdocvar{P} \coqdocvar{Q}) $\rightarrow$ \coqdoceol
\coqdocindent{2.00em}(\coqdocvar{b\_product} \coqdocvar{R} \coqdocvar{S}) $\equiv$ (InteriorG Z).
\end{coqdoccode}
\end{theorem}
It follows that the geometric space with balls can be interpreted as a topology $\langle balls, \mathbb{T}_{MG}\rangle$, based on the \coq{InteriorG} definition. To complete $\mathbb{T}_{MG}$, we first add the boundary definition already suggested by \cite{Jaskowski1948} and implicit in Tarski's definition of equidistance \cite{Tarski1956b}. Informally, it says that a point \(P\) belongs to the boundary of a given solid \(Q\) iff \(P\) is both not a part of \(Q\) and not external to \(Q\):
\begin{definition}[\coq{boundary}]\label{boundary}\mbox{}\\
$\forall \: P \: Q \: : N, \: \eta \: P \: (boundary\: Q) \: \triangleq \: (\eta \: Q \: Solid \; \wedge \; \eta \: P \: balls \; \wedge \; \forall \: Y, \; \eta \: Y \: balls \; \wedge \; \eta \: P \: (\coq{Point}\,Y) \rightarrow$ \\ $\neg \eta \: Y \: (pt \: Q) \; \wedge \;
\neg \eta \: Y \: (ext \: Q)).$
\end{definition}
Besides the boundary definition, a related mereological definition was already proposed in \cite{Lejewski1954} and known as $overlap$:
 \begin{definition}[\coq{overlap}]\label{overlap}\mbox{}
$\forall \; P \; Q, \; \eta \: P \: (\mathtt{overlap} \: Q) \: \triangleq \: \mathtt{Individual} \: (P) \; \wedge \; \exists \: R, (\eta \: R \; (pt \: P) \: \wedge \: \eta \: R \; (pt \: Q)).$
\end{definition}
Then, it is routine to prove that the overlap definition coincides with the negation of $ext$:
\begin{lemma}[\coq{ov\_not\_ext}]\label{ovNotExt}\mbox{}
\begin{coqdoccode}
\ensuremath{\forall} \coqdocvar{P} \coqdocvar{Q}, \coqdocdefinition{Individual} \coqdocvar{P} \ensuremath{\land} \coqdocdefinition{Individual} \coqdocvar{Q} \ensuremath{\rightarrow} \coqdocdefinition{$\eta$} \coqdocvar{P} (\coqdocdefinition{overlap} \coqdocvar{Q}) $\leftrightarrow$ $\neg$ \coqdocdefinition{$\eta$} \coqdocvar{P} (\coqdocdefinition{ext} \coqdocvar{Q})).
\end{coqdoccode}
\end{lemma}
It results that the last term in definition \ref{boundary} can be substituted by a positive term, i.e., $\eta \: Y \: (\mathtt{overlap} \: Q)$\footnote{This fact has also been noticed in \cite{Gruszczynski2008}.}.
Using the boundary (denoted \(\delta\)), we can prove among others that, (i) the boundary of an empty solid does not exist, (ii) the space has no boundary, (iii) a solid and its mereological complement share the same boundary, (iv) interior points of a solid and boundary points are disjoint and (v) boundary is not a solid.
\begin{lemma}[\coq{no\_empty\_bound}]\label{noEmptyBound}\mbox{}
\begin{coqdoccode}
($\delta$ \coqdocvar{$\varLambda$}) $\approx$ \coqdocvar{$\varLambda$}.
\end{coqdoccode}
\end{lemma}
\begin{lemma}[\coq{no\_Gspace\_bound}]\label{noGspaceBound}\mbox{}
\begin{coqdoccode}
\ensuremath{\forall} \coqdocvar{P}, \coqdocdefinition{$\eta$} \coqdocvar{P} ($\delta$ \coqdocvar{Gspace}) \ensuremath{\rightarrow} \coqdocdefinition{$\eta$} \coqdocvar{P} \coqdocvar{$\varLambda$}.
\end{coqdoccode}
\end{lemma}
\begin{lemma}[\coq{compl\_bound}]\label{complbound}\mbox{}
\begin{coqdoccode}
\ensuremath{\forall} \coqdocvar{P} \coqdocvar{Q}, \coqdocdefinition{$\eta$} \coqdocvar{P} (\coqdocdefinition{relCompl} \coqdocvar{Q} \coqdocvar{Gspace}) \ensuremath{\rightarrow} $\delta$ \coqdocvar{Q} $\approx$ $\delta$ \coqdocvar{P}.
\end{coqdoccode}
\end{lemma}
\begin{lemma}[\coq{bound\_int\_disj}]\label{boundIntDisj}\mbox{}
\begin{coqdoccode}
\ensuremath{\forall} \coqdocvar{P} \coqdocvar{Q}, \coqdocdefinition{$\eta$} \coqdocvar{P} (\coqdocdefinition{sat\_InteriorPoint} \coqdocvar{Q}) \ensuremath{\land} \coqdocdefinition{$\eta$} \coqdocvar{P} ($\delta$ \coqdocvar{Q}) \ensuremath{\rightarrow} \coqdocdefinition{$\eta$} \coqdocvar{P} \coqdocvar{$\varLambda$}.\coqdoceol
\end{coqdoccode}
\end{lemma}
\begin{lemma}[\coq{bound\_not\_sol}]\label{notSol}\mbox{}
\begin{coqdoccode}
\ensuremath{\forall} \coqdocvar{P} \coqdocvar{Q}, $\eta$ \coqdocvar{Q} \coqdocdefinition{Solid} \ensuremath{\land} $\eta$ \coqdocvar{P} ($\delta$ \coqdocvar{Q}) \ensuremath{\rightarrow} $\neg$(($\delta$ \coqdocvar{Q}) $\approx$ \coqdocvar{Solid}).\coqdoceol
\end{coqdoccode}
\end{lemma}
The mereological complement satisfies the definition involved in \ref{Btopo} in which \(R\) is \(Gspace\). Since solids coincide with open sets, the last lemma proves that the boundary is not open which prevents it to appear in terms in which an open entity is required. In summary, $\mathbb{T}_{MG}$ is a topology of open solids and complies with many situations described, in mereotopology.
\subsection{An application example: RCC8}
As an application example, we are able to express all relations involved in the Region Connection Calculus (RCC8) \cite{Randell1992} and to reason with them. However, it is well known that in the RCC theory it is not possible to distinguish between regions that are open, closed or neither but have the same closure. The purpose here is to express these binary relations with mereological primitives (e.g., $pt$, $ppt$, $ext$) and the boundary relation as a substitute to the "connect" relation on which RCC8 is based. The following examples show how RCC8-style relations can be expressed as definitions on open solids. These eight relations respectively include $EQ$, $DC$ (disconnected), $EC$ (externally connected), $PO$ (partial overlap), $TPP$ (Internal tangential proper part), $NTPP$ (Internal non-tangential proper part), $TPPi$ (Internal tangential proper part inverted) and $NTPPi$ (Internal non-tangential proper part inverted).
\begin{definition}[\coq{EQ}]\label{EQ}\mbox{}
$\forall \; P \; Q, \eta \: P \: (\mathtt{EQ} \: Q) \: \triangleq \; P \; \equiv \; Q.$
\end{definition}
\begin{definition}[\coq{DC}]\label{DC}\mbox{}
$\forall \; P \; Q, \; \eta \: P \: (\mathtt{DC} \: Q) \: \triangleq \: \eta \: P \; (ext \: Q)
\wedge \; \neg \exists \: R, (\eta \: R \; (\delta \: P) \: \wedge \: \eta \: R \; (\delta \: Q)).$
\end{definition}
\begin{definition}[\coq{EC}]\label{EC}\mbox{}
$\forall \; P \; Q, \; \eta \: P \: (\mathtt{EC} \: Q) \: \triangleq \; \eta \: P \; (ext \: Q) \; \wedge \; \exists \: R, (\eta \: R \; (\delta \: P) \: \wedge \: \eta \: R \; (\delta \: Q)).$ 
\end{definition}
\begin{definition}[\coq{PO}]\label{PO}\mbox{}
$\forall \; P \; Q, \; \eta \: P \: (\mathtt{PO} \: Q) \: \triangleq \: \neg \eta \: P \; (pt \: Q) \: \wedge \: \neg \eta \: Q \; (pt \: P) \: \wedge \; \eta \: P \: (\mathtt{overlap} \: Q).$ 
\end{definition}
\begin{definition}[\coq{TPP}]\label{TPP}\mbox{}
$\forall \; P \; Q, \; \eta \: P \: (\mathtt{TPP} \: Q) \: \triangleq \: \eta \: P \; (ppt \: Q) \; \wedge \: \exists \: R, (\eta \: R \; (\delta \: P) \: \wedge \: \eta \: R \; (\delta \: Q)).$ 
\end{definition}
\begin{definition}[\coq{NTPP}]\label{NTPP}\mbox{}
$\forall \; P \; Q, \; \eta \: P \: (\mathtt{NTPP} \: Q) \: \triangleq \: \eta \: P \; (ppt \: Q) \; \wedge \: \neg \exists \: R, (\eta \: R \; (\delta \: P) \: \wedge \: \eta \: R \; (\delta \: Q)).$ 
\end{definition}
\begin{definition}[\coq{TPPi}]\label{TPPi}\mbox{}
$\forall \; P \; Q, \; \eta \: P \: (\mathtt{TPPi} \: Q) \: \triangleq \: \eta \: Q \; (TPP \: P).$ 
\end{definition}
\begin{definition}[\coq{NTPPi}]\label{NTPPi}\mbox{}
$\forall \; P \; Q, \; \eta \: P \: (\mathtt{NTPPi} \: Q) \: \triangleq \: \eta \: Q \; (NTPP \: P).$ 
\end{definition}
All these relations are pairwise disjoint and enhance the precision in RCC8 since they are not primitive but rather rely on the formalism introduced here. Furthermore, using the list of definitions it becomes possible to prove the so-called elements in the consistency-based composition table (see e.g., \cite{Sanjiang2003}). First, some rewriting lemmas are required to ease further theorems, e.g.:
\begin{lemma}[\coq{TPP\_to\_NTPP}]\label{TPPnotNTPP}\mbox{}
\begin{coqdoccode}
\ensuremath{\forall} \coqdocvar{P} \coqdocvar{Q}, \coqdocdefinition{$\eta$} \coqdocvar{P} (\coqdocdefinition{TPP} \coqdocvar{Q}) $\leftrightarrow$ \coqdocdefinition{$\eta$}  \coqdocvar{P}  (\coqdocdefinition{ppt} \coqdocvar{Q}) \ensuremath{\land} $\neg$ $\eta$ \coqdocvar{P} (\coqdocdefinition{NTPP} \coqdocvar{Q}).\coqdoceol
\end{coqdoccode}
\end{lemma}
Then, cell entries of the composition table are proved using definitions, rewriting lemmas and mereological theorems of $\lambda$-MM. Let us consider for example the cell entry $TPP$ to $PO$ in the composition table. It corresponds to the expression 
$$\forall \; X \; Y \; Z, \; \: TPP \; (X, Y) \: \wedge PO \; (Y, Z) \rightarrow \varphi \; (X, Z)$$
where $\varphi$ is a list of possible and exhaustive relations which, in that case, corresponds to $\lbrace \; \mathtt{DC}, \; \mathtt{EC}, \; \mathtt{PO}, \; \mathtt{TPP}, \; \mathtt{NTPP}\rbrace$.
The theorem is written as follows:
\begin{theorem}[\coq{TPP\_to\_PO}]\label{TPPtoPO}\mbox{}
\begin{coqdoccode}
 \ensuremath{\forall} \coqdocvar{P} \coqdocvar{Q} \coqdocvar{R}, \coqdocvar{$\eta$} \coqdocvar{P} \coqdocvar{Solid} \ensuremath{\land} \coqdocvar{$\eta$} \coqdocvar{Q} \coqdocvar{Solid} \ensuremath{\land} \coqdocvar{$\eta$} \coqdocvar{R} \coqdocvar{Solid} \ensuremath{\rightarrow} \coqdocvar{$\eta$} \coqdocvar{P} (\coqdocdefinition{TPP} \coqdocvar{Q}) \ensuremath{\rightarrow} \coqdoceol
\coqdocindent{2.00em} \coqdocvar{$\eta$} \coqdocvar{Q} (\coqdocdefinition{PO} \coqdocvar{R}) \ensuremath{\rightarrow} \coqdocvar{$\eta$} \coqdocvar{P} (\coqdocdefinition{DC} \coqdocvar{R}) \ensuremath{\lor} \coqdocvar{$\eta$} \coqdocvar{P} (\coqdocdefinition{EC} \coqdocvar{R}) \ensuremath{\lor} \coqdocvar{$\eta$} \coqdocvar{P} (\coqdocdefinition{PO} \coqdocvar{R}) \ensuremath{\lor}  \coqdocvar{$\eta$} \coqdocvar{P} (\coqdocdefinition{TPP} \coqdocvar{R}) \ensuremath{\lor} \coqdocvar{$\eta$} \coqdocvar{P} (\coqdocdefinition{NTPP} \coqdocvar{R}).\coqdoceol
\end{coqdoccode}
\end{theorem}
\paragraph{Proof sketch.}
The proof proceeds by a finite case analysis on the relation between \(P\)
and \(R\). From \(\mathtt{TPP}(P,Q)\) and \(\mathtt{PO}(Q,R)\), a
preliminary lemma first shows that \(R\) is not a part of \(P\). If \(P\)
is external to \(R\), the definitions of \(\mathtt{DC}\) and
\(\mathtt{EC}\) reduce the conclusion to whether the boundaries of \(P\)
and \(R\) meet. If \(P\) is not external to \(R\), the equivalence between
non-externality and overlap gives \(\mathtt{overlap}(P,R)\). A second case
distinction asks whether \(P\) is a part of \(R\). If it is not, the
definition of \(\mathtt{PO}\) applies. If it is, the previous fact that
\(R\) is not a part of \(P\) turns parthood into proper parthood; a final
split on boundary intersection then yields either \(\mathtt{TPP}(P,R)\)
or \(\mathtt{NTPP}(P,R)\). Thus the five possible RCC8 outcomes are
exhaustive.
\subsection{Formal geometric background}
\label{subsec:formal-geometric-background}

The Coq development separates the topological construction from the
geometric assumptions used to model Tarski's geometry of solids. Besides
the underlying Le{\'s}niewskian nominal and mereological primitives, the
explicitly geometric part of the library assumes the primitive predicate
\texttt{isBall}, the non-emptiness of the ball domain, the compatibility
between concentricity and point-class representation, closure of
solids under parts, the formal axiom \texttt{Three\_points}, and an
extensionality principle for common surface behaviour.
These assumptions form the geometric background of what will be called
the Tarskian regional construction. The point of the next sections is to
make explicit the representative domain induced from this regional
framework, so that topological operations can be specified on the
appropriate point-class objects rather than on regional solids
themselves.
The axiom \texttt{Three\_points} is included here because it belongs to
the regional geometric background of the library: it states the existence
of a \texttt{btw} representative for non-concentric balls. The present
paper does not develop the theory of betweenness itself; \texttt{btw}
serves only as part of the inherited Tarskian geometry of solids.

\begin{center}
\footnotesize
\setlength{\tabcolsep}{4pt}
\begin{tabular}{@{}p{0.36\textwidth}p{0.32\textwidth}p{0.22\textwidth}@{}}
\hline
Coq assumption & Informal role & Use \\
\hline
\texttt{isBall} & primitive ball predicate & geometric signature \\
\texttt{ball\_exists} & non-empty ball domain & representative existence \\
\texttt{Concent\_same\_tarski\_point} &
concentricity preserves point-class identity & point reconstruction \\
\texttt{part\_of\_solid} & parts of solids are solids & regional mereology \\
\texttt{Three\_points} & existence for \texttt{btw} on non-concentric balls & regional geometry \\
\texttt{same\_center\_common\_surface} &
surface/extensionality principle & representative geometry \\
\hline
\end{tabular}
\end{center}

\section{From Regional Opens to the Point-Class Domain}
\label{sec:regional-point-domain}

The regional topology developed above provides the geometric and
mereological setting of the construction. It describes solids, their
parts, and the regional operations inherited from Tarski's geometry. The
new construction has a different carrier: it concerns reconstructed
locations represented by concentric families of balls. The task of the
following sections is therefore to construct a topology on these
representatives and then its neighbourhood closure.

\paragraph{Why the quotient is unavoidable.}
In Tarski's geometry of solids, a reconstructed point is not an individual
ball. It is the entire family of balls with one center. No member of this
family has privileged ontological status. Consequently, every topological
notion intended to concern reconstructed points must be invariant under
replacement of a ball by a concentric representative. Otherwise,
neighbourhoods or closure would depend on the accidental choice of a radius
rather than on the reconstructed location.

The quotient induced by \SameCenter{} is thus the canonical carrier
determined by the reconstruction. Conversely, a topology defined directly
on individual ball representatives would distinguish objects that represent
exactly the same reconstructed point. The formal development represents
this quotient setoid-theoretically: ball representatives remain in the
original nominal domain, while equality, membership, and topological
predicates are required to respect \SameCenter{}. This preserves both
Tarski's representative-based reconstruction and the single nominal
ontology inherited from \(\lambda\)-MM.

\paragraph{Why setoids rather than quotient types?}
The quotient is conceptual, but it is not reified as a new Coq type. This
choice is not merely an implementation convenience. Tarski's
reconstruction does not introduce points as a new ontological category:
the point-like object is given extensionally by a family of concentric
balls. A setoid retains the ball representatives in the nominal domain and
identifies them only for the purposes relevant to equality, membership,
neighbourhood, and closure. Introducing a separate quotient type would be
mathematically possible, but it would less directly reflect the
representative-based and non-reifying orientation of the present
formalization.

\paragraph{From regional generators to a point-class carrier.}
The regional and point-class registers must therefore be kept distinct.
\coq{Gspace} is the ambient regional individual generated by balls, whose
elements are treated mereologically. The point-class carrier is not
\coq{Gspace} itself, but the representative plural \coq{balls} modulo
\SameCenter{}. A regional ball \(Q\) will provide a basic plural of
point-class representatives; it is not itself an open set of the new
topology. The next section makes this passage explicit, first by defining
the setoid and extensional point-class plurals, then by showing that the
ball-generated plurals satisfy the two axioms of a topological basis, and
only afterwards by defining the open plurals and their closure.

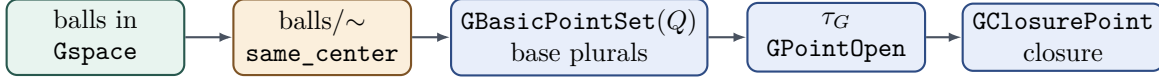
\begin{figure}[t]
\centering
\begin{tikzpicture}[
  block/.style={rectangle, rounded corners, draw, thick, align=center,
                minimum width=2.35cm, minimum height=.75cm, font=\small},
  regional/.style={block, draw=RSLRegionalInk, fill=RSLRegionalFill},
  point/.style={block, draw=RSLPointInk, fill=RSLPointFill},
  bridge/.style={block, draw=RSLBridgeInk, fill=RSLBridgeFill},
  flow/.style={-{Latex[length=1.7mm]}, thick, draw=RSLArrow}
]
\node[regional] (balls) at (0,0) {balls in\\\(\coq{Gspace}\)};
\node[bridge] (quotient) at (3.0,0) {\(\mathrm{balls}/{\sim}\)\\\(\SameCenter\)};
\node[point] (basic) at (6.4,0) {\(\GBasicPointSetName(Q)\)\\base plurals};
\node[point] (topology) at (9.8,0) {\(\tau_G\)\\\(\GPointOpenName\)};
\node[point] (closure) at (12.8,0) {\(\GClosurePointName\)\\closure};
\draw[flow] (balls) -- (quotient);
\draw[flow] (quotient) -- (basic);
\draw[flow] (basic) -- (topology);
\draw[flow] (topology) -- (closure);
\end{tikzpicture}
\caption{The verified progression from regional ball generators to the
ball-generated point-class topology and its neighbourhood closure.}
\label{fig:geometric-topology}
\end{figure}

\section{The Ball-Generated Topology of Reconstructed Points}
\label{sec:geometric-point-topology}

\subsection{Reconstructed points as a setoid}
\label{subsec:geometric-setoid}

Let \(\coq{Gspace}\) denote the ambient mereogeometric region. Every ball is a part of it:
\begin{lemma}[\coq{in\_Gspace}]\label{lem:new-ing-space}
\[
\forall A,\qquad
\eta A\,\coq{balls}\ \longrightarrow\ \eta A\,(\coq{pt}\,\coq{Gspace}).
\]
\end{lemma}
Thus \(\coq{Gspace}\) is the regional ambient individual, whereas the
carrier for reconstructed points is the plural of ball representatives.
\begin{definition}[\SameCenter]\label{def:new-samecenter}
\[
\forall P\,Q:N,\qquad
\SameCenter(P,Q)
\quad\triangleq\quad
\coq{Point}(P)\approx\coq{Point}(Q).
\]
\end{definition}
For balls \(P,Q\), write \(P\sim Q\) for \(\SameCenter(P,Q)\). The Coq
lemmas \coq{same\_center\_refl}, \coq{same\_center\_sym}, and
\coq{same\_center\_trans} establish that this is an equivalence relation.
A reconstructed point is therefore a \emph{concentric point-class}
\[
[P]_\sim=\{\,Q\mid \eta Q\,\coq{balls}\ \wedge\ Q\sim P\,\}.
\]
The reconstructed-point domain is presented extensionally as
\[
X_G=\coq{balls}/{\sim}.
\]
No new Coq type is introduced: subsets of \(X_G\) are nominal plurals of
ball representatives whose membership is stable under \SameCenter{}.

\subsection{Point-class plurals and extensional equality}
\label{subsec:geometric-pointset-plurals}

A plural of reconstructed points is represented by a nominal name whose
members are ball representatives. Since a point-class has many
representatives, equality and inclusion are tested extensionally only on balls:
\begin{definition}[\PointSetEq]\label{def:new-pointseteq}
\[
\PointSetEq(u,v)
\quad\triangleq\quad
\forall B,\ \eta B\,\coq{balls}\rightarrow
\bigl(\eta B\,u\leftrightarrow\eta B\,v\bigr).
\]
\end{definition}
\begin{definition}[\PointSetIncl]\label{def:new-pointsetincl}
\[
\PointSetIncl(u,v)
\quad\triangleq\quad
\forall B,\ \eta B\,\coq{balls}\rightarrow
\eta B\,u\rightarrow\eta B\,v.
\]
\end{definition}
\begin{definition}[\PointSetUnion]\label{def:new-pointsetunion}
\[
\eta B\,(u\cup_{pt}v)
\quad\triangleq\quad
\eta B\,\coq{balls}\wedge\bigl(\eta B\,u\vee\eta B\,v\bigr).
\]
\end{definition}
\begin{definition}[\PointSetStable]\label{def:new-pointsetstable}
\[
\PointSetStable(u)
\quad\triangleq\quad
\forall B\,C,\ \SameCenter(B,C)\rightarrow
\eta B\,u\rightarrow\eta C\,u.
\]
\end{definition}
\begin{definition}[\PointSetMeet]\label{def:new-pointsetmeet}
\[
\PointSetMeet(u \, v)
\quad\triangleq\quad
\exists B\,C,\ \eta B\,\coq{balls} \, \wedge \, \eta C\,\coq{balls} \, \wedge \,\SameCenter(B,C) \wedge
\eta B\,u \, \wedge \, \eta C\,v.
\]
\end{definition}
We write \(=_{pt}\), \(\subseteq_{pt}\), and \(\cup_{pt}\) for point-class equality, inclusion, and union, respectively. Stability ensures that a plural depends on a concentric point-class rather than on the radius of a chosen representative.
The last definition \(\PointSetMeet(u \, v)\) corresponds to the existence of two concentric balls, one in u, the other in v.
\subsection{Ball-generated topological base}
\label{subsec:geometric-basic}

We first introduce the interior point \(\PointClassInteriorName\) of a set of points \(q\):
\begin{definition}[\PointClassInteriorName]\label{def:ginterior}
\[
\begin{aligned}
\eta P\,\PointClassInteriorName(q)
\quad\triangleq\quad&
\eta P\,\coq{balls}\ \wedge\ \exists C\,(\eta C\,\coq{balls}\wedge
\SameCenter(P,C)\wedge\eta C\,q).
\end{aligned}
\]
\end{definition}

For each ball \(Q\), the formalization defines the basic point-plural
\(\GBasicPointSetName(Q)\). Its membership condition is:
\begin{definition}[\GBasicPointSetName]\label{def:gbasic}
\[
\begin{aligned}
\eta P\,\GBasicPointSetName(Q)
\quad\triangleq\quad&
\eta Q\,\coq{balls}\ \wedge\ \eta P\,\coq{balls}\\
&{}\wedge\ \exists C\,\bigl(\eta C\,\coq{balls}\wedge
\SameCenter(P,C)\wedge\eta C\,(\coq{pt}\,Q)\bigr).
\end{aligned}
\]
\end{definition}

Thus \(\GBasicPointSetName(Q)\) contains precisely the representatives of point-classes having a ball representative regionally contained in \(Q\). It is routine to see that, provided that \(Q\) is a ball, \(\GBasicPointSetName(Q)\) coincides with the interior parts of \(Q\):
\begin{lemma}[\coq{GBasicPointSet\_as\_interior}]\label{basToInt}
\[
\forall P Q,\quad
\eta Q\,\coq{balls}\, \longrightarrow\
(\eta P\,\GBasicPointSetName (Q) \longleftrightarrow\ \eta\, P\, \PointClassInteriorName (\coq{pt}\, Q)).
\]
\end{lemma}
We first have to prove that \(\GBasicPointSetName(Q)\) behave like a topological base \(\mathcal{B}\). For that purpose, \(\GBasicPointSetName(Q)\) must verify two intrinsic conditions (often called the intrinsic basis axioms): (i) a cover axiom, 
asserting that for every point \(x \in X\), there is at least one basic set \(B \in \mathcal{B}\) such that \(x \in B\) and (ii) an intersection axiom, which claim that for any two basic sets \(B_1, B_2 \in \mathcal{B}\) and any point \(x \in B_1 \cap B_2\), there must exist a third basic set \(B_3 \in \mathcal{B}\) such that \(x \in B_3 \subseteq (B_1 \cap B_2)\).

For readability, we write \(\mathcal{G}(Q)=\GBasicPointSetName(Q)\). The family of these plurals, denoted \(\mathcal{G}\), is the base from which the geometric topology is generated. The formalization proves the usual cover and intersection conditions.  

\begin{theorem}[\coq{GBasicPointSet\_is\_topological\_base}]
\label{thm:gbasic-topological-base}
\[
\forall P,\quad
\eta P\,\coq{balls}\ \longrightarrow\
\eta P\,\mathcal{G}(P),
\]
and
\[
\begin{aligned}
\eta P\,\mathcal{G}(B)\wedge
\eta P\,\mathcal{G}(Q)
\ \longrightarrow\ 
\exists D,\ &\eta D\,\coq{balls}\wedge
\eta P\,\mathcal{G}(D)\\
&{}\wedge\PointSetIncl(\mathcal{G}(D),\mathcal{G}(B))\\
&{}\wedge\PointSetIncl(\mathcal{G}(D),\mathcal{G}(Q)).
\end{aligned}
\]
\end{theorem}
The first clause says that every ball representative is covered by a basic plural. The second says that whenever a representative belongs to two
basic plurals, it belongs to a third one included in both.  Hence
\(\{\GBasicPointSetName(Q)\mid\eta Q\,\coq{balls}\}\) is a topological
basis on the setoid carrier \(X_G\).  The following definition presents
the opens generated by this basis directly in the nominal language.

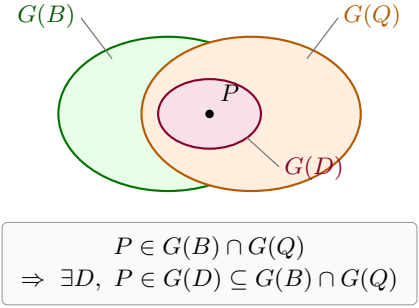
\begin{figure}[htbp]
\centering
\begin{tikzpicture}[
  x=1cm,
  y=1cm,
  >=Latex,
  font=\small,
  regionball/.style={draw=green!45!black, fill=green!10, line width=0.8pt},
  witnessball/.style={draw=blue!65!black, fill=blue!16, line width=0.8pt},
  repball/.style={draw=blue!70!black, dashed, line width=0.9pt},
  basic/.style={draw=blue!65!black, fill=blue!10, line width=0.8pt},
  basicb/.style={draw=green!45!black, fill=green!10, line width=0.8pt},
  basicq/.style={draw=orange!70!black, fill=orange!13, line width=0.8pt},
  refined/.style={draw=purple!70!black, fill=purple!12, line width=0.8pt},
  callout/.style={draw=black!40, fill=black!2, rounded corners=2pt, align=center, inner xsep=7pt, inner ysep=5pt, font=\footnotesize},
  labelnode/.style={font=\footnotesize, align=center, inner sep=1pt},
  leader/.style={draw=black!55, line width=0.45pt}
]

\begin{scope}
  \node[font=\bfseries] at (0.3,2.48) {(a) Membership in a basic neighbourhood};

  \draw[regionball] (0,0) circle (1.35);
  \draw[leader] (-1.63,1.28) -- (-1.05,0.84);
  \node[labelnode, text=green!35!black, anchor=east] at (-1.66,1.30) {$Q$};

  \coordinate (pcenter) at (0.55,0.15);
  \draw[witnessball] (pcenter) circle (0.38);
  \draw[leader] (1.16,-0.38) -- (0.85,-0.08);
  \node[labelnode, text=blue!65!black, anchor=west] at (1.19,-0.40) {$C$};
  \node[labelnode, text=blue!65!black, anchor=west] at (1.55,-0.02) {$C \subseteq Q$};

  \draw[repball] (pcenter) circle (1.05);
  \draw[leader] (1.94,1.28) -- (1.33,0.90);
  \node[labelnode, text=blue!70!black, anchor=west] at (1.97,1.30) {$P$};
  \node[labelnode, text=blue!70!black, anchor=west] at (1.74,0.48) {$P \nsubseteq Q$};

  \fill[black] (pcenter) circle (1.4pt);
  \draw[leader] (0.08,-0.60) -- (0.48,0.06);
  \node[labelnode, anchor=east] at (0.03,-0.62) {common center};

  \node[callout] at (0.25,-1.98)
    {$C \sim P$ and $C \subseteq Q$\\[1pt]
     hence $\eta\,P\,(\mathrm{GBasicPointSet}(Q))$};
\end{scope}

\begin{scope}[xshift=6.80cm]
  \node[font=\bfseries] at (0,2.48) {(b) Intersection};
  \node[font=\scriptsize, text=black!65] at (0,2.08)
    {$G(R)=\mathrm{GBasicPointSet}(R)$};

  \draw[basicb] (-0.55,0) ellipse (1.45 and 1.02);
  \draw[basicq] (0.55,0) ellipse (1.45 and 1.02);
  \draw[refined] (0,0) ellipse (0.68 and 0.46);

  \draw[leader] (-1.70,1.27) -- (-1.29,0.73);
  \node[labelnode, text=green!35!black, anchor=east] at (-1.73,1.29) {$G(B)$};
  \draw[leader] (1.70,1.27) -- (1.29,0.73);
  \node[labelnode, text=orange!70!black, anchor=west] at (1.73,1.29) {$G(Q)$};
  \draw[leader] (0.92,-0.70) -- (0.49,-0.31);
  \node[labelnode, text=purple!70!black, anchor=west] at (0.95,-0.72) {$G(D)$};

  \fill[black] (0,0) circle (1.6pt);
  \node[labelnode, anchor=south west] at (0.10,0.13) {$P$};

  \node[callout] at (0,-1.98)
    {$P \in G(B) \cap G(Q)$\\[1pt]
     $\Rightarrow\ \exists D,\ P \in G(D) \subseteq G(B)\cap G(Q)$};
\end{scope}

\end{tikzpicture}
\caption{Ball-generated basic point-plurals. Panel (a) illustrates membership in \(\GBasicPointSetName(Q)\) through a concentric witness ball contained in \(Q\). Panel (b) illustrates the local refinement of two intersecting basic point-plurals.}
\label{fig:gbasic-pointset-basis}
\end{figure}

\subsection{Geometric point-open plurals}
\label{subsec:geometric-open}

A plural is geometrically open\footnote{Throughout Sections~7--8, the qualifier ``geometric'' indicates that a notion is generated from the regional ball structure of Tarski's geometry of solids, in contrast with the individual-name topology of Section~4; it does not refer to geometric topology in the sense of manifold theory.} exactly when it is stable under
\SameCenter{}, contains only balls, and each of its members has a
ball-generated basic plural contained in it:
\begin{definition}[\GPointOpenName]\label{def:gpointopen}
\[
\begin{aligned}
\GPointOpenName(u)\quad\triangleq\quad&
\PointSetStable(u)\ \wedge\
\bigl(\forall P,\eta P\,u\rightarrow\eta P\,\coq{balls}\bigr)\\
&{}\wedge\
\forall P,\eta P\,u\rightarrow\exists Q,\eta Q\,\coq{balls}\\
&{}\wedge\eta P\,\GBasicPointSetName(Q)
\wedge\PointSetIncl(\GBasicPointSetName(Q),u).
\end{aligned}
\]
\end{definition}
The ball requirement on the witness \(Q\) is essential. If arbitrary
saturated point-class plurals were allowed as witnesses, every stable
plural could witness itself and the construction would collapse to a
discrete topology.

Each element in \(\mathcal{G}\) is geometrically open:
\[
\eta Q\,\coq{balls}\quad\longrightarrow\quad
\GPointOpenName(\GBasicPointSetName(Q)),
\]
which is the Coq lemma \coq{GBasicPointSet\_is\_GPointOpen}.  Consequently the displayed definition is exactly the standard base-generated topology \(\langle \tau_G, \mathcal{G} \rangle\) expressed through stable nominal plurals. 

\subsection{The topology theorem}
\label{subsec:geometric-topology-theorem}

Let
\[
\tau_G=\{\,u:N\mid \GPointOpenName(u)\,\}.
\]
The formal development proves independently that \(\tau_G\) contains the
empty plural and the full representative domain \(\coq{balls}\), and is
closed under binary intersections and arbitrary indexed unions.
\begin{theorem}[\coq{GPointOpen\_is\_topology}]
\label{thm:gpointopen-topology}
\[
\begin{aligned}
&\GPointOpenName(\Lambda),\\
&\GPointOpenName(\coq{balls}),\\
&\GPointOpenName(u)\wedge\GPointOpenName(v)
\ \longrightarrow\ \GPointOpenName(u\cap v),
\end{aligned}
\]
and, for every metatheoretic family \(U:I\rightarrow N\),
\[
\bigl(\forall i:I,\GPointOpenName(U_i)\bigr)
\ \longrightarrow\
\GPointOpenName\Bigl(\bigcup_{i:I}U_i\Bigr).
\]
\end{theorem}
Hence \(\tau_G\) is a topology on the setoid presentation \(X_G\). The
indexed union is metatheoretic and introduces no powerset object into the Le{\'s}niewskian ontology.

\subsection{Conditional non-discreteness}
\label{subsec:geometric-nondiscrete}

Non-discreteness requires an additional geometric richness hypothesis. The definition \(\BallCenterBundleName \) states the local geometric condition needed below: every regional ball contains another ball with a different center, i.e., no elementary regional neighbourhood is reduced to a single center.
\begin{definition}[\BallCenterBundleName]\label{def:ball-richness}
\[
\begin{aligned}
\coq{BallCenterBundle}
\quad\triangleq\quad&
\forall B,\eta B\,\coq{balls}\rightarrow
\exists C,\eta C\,\coq{balls}\wedge\eta C\,(\coq{pt}\,B)
\wedge\neg\SameCenter(C,B).
\end{aligned}
\]
\end{definition}
Let us consider \(\GPointSingleName \; P\) as the plural name which stands for all balls falling into the same reconstructed point class generated by \( P\).

\[
\coq{GPointSingle} (P)
\quad\triangleq\quad
\PointClassInteriorName\, (\coq{Point}\, P).
\]
It is easy to prove that all balls in \(\GPointSingleName \) have the same center:
\begin{lemma}[\coq{is\_GPointSingle}]
\label{gpointsinglecenter}
\[
\forall C \, P, \quad\ \eta P\,\coq{balls}\rightarrow
\eta \, C\,\coq{GPointSingle} (P)
\quad \longleftrightarrow \quad
 \eta C\,\coq{balls} \wedge\, \SameCenter\, C \, P .
\]
\end{lemma}
Single point class discreteness is implied by the literal definition above, once they are shown to be geometric point-plurals. A geometric point-plural is a stable set of rebuilt points whose members all belong to the domain [balls]:
\[
\coq{GPointPlural} (u)
\quad\triangleq\quad
\PointSetStable\, u \wedge (\forall P, \eta P\, u \rightarrow \eta P\,\coq{balls}).
\]
The topology is discrete if for all plural \(u\), a proof of \(\coq{GPointPlural} \; u\) yields a proof for \(\GPointOpenName \; u\). While the left member can be easily proved:
\begin{lemma}[\coq{GPointSingle\_is\_GPointPlural}]
\label{gpointsingleplural}
\[
\forall P, \quad\ \coq{GPointPlural} (\coq{GPointSingle} (P)).
\]
\end{lemma}
the right member cannot.

Now, suppose that \(\GPointSingleName \; P\) is open. Since \(P\) belongs to its single class, openness generates a base neighbourhood \(\GBasicPointSetName(Q)\) including \(P\) and enclosed inside \(\GPointSingleName \; P\). However, any membership to this base neighbourhood means that there exists a ball \(D\) concentric with \( P\), inside \(Q\). According to definition \(\BallCenterBundleName \), \(D\) includes a sub-ball \(C\) having a distinct center. Since \(C \subset D \subset Q\), the point generated by \(C\) also belongs to the base \(\GBasicPointSetName(Q)\). Since such a base neighbourhood is supposed included in \(\GPointSingleName \; P\), it results that \(C\) has the same center as \(P\), and then, like \(D\), which leads to a contradiction.

The corresponding Coq theorem is \coq{BallCenterBundle\_no\_open\_point\_singleton}.
\begin{theorem}[Non-discreteness from local richness]
\label{thm:gpoint-topology-not-discrete}
\[
\coq{BallCenterBundle} \longrightarrow\
\forall P,\quad
\eta P\,\coq{balls}\ \longrightarrow\ 
\neg\coq{GPointOpen} (\coq{GPointSingle} (P)).
\]
\end{theorem}
The theorem \coq{BallCenterBundle\_implies\_GPointTopology\_not\_discrete} then establishes that the topology is not discrete.
\subsection{Hausdorff separation}
\label{subsec:geometric-hausdorff}

For distinct concentric point-classes, Tarski's \coq{Three\_points} axiom together with the defining characterization of betweenness yield exterior ball
representatives. The associated basic point-plurals are disjoint. Proofs rely on the \(\coq{ext}\) operator of $\lambda$-MM. The first step provides a proof that \(\coq{GBasicPointSet}\) of external generated points are disjoints:
\begin{lemma}[\coq{exterior\_GBasicPointSet\_disjoint}]
\label{gpointdisjoint}
\[
\begin{aligned}
\forall B \, C,  \eta \, B\,\coq{balls}\ \rightarrow\  &\eta \, C\,\coq{balls}\ \rightarrow\ \eta \, B\, (\coq{ext} C) \rightarrow\ \\ &\neg\PointSetMeet (\GBasicPointSetName B) (\GBasicPointSetName C).
\end{aligned}
\]
\end{lemma}
The second step shows that distinct balls generate corresponding points plurals in which some balls are external:
\begin{lemma}[\coq{distinct\_centers\_have\_exterior\_basic\_balls}]
\label{extbasepoints}
\[
\begin{aligned}
\forall P \, Q,  \eta \, P\,\coq{balls}\ \rightarrow\  &\eta \, Q\,\coq{balls}\ \rightarrow\ \neg\SameCenter\, P \, Q \rightarrow\ \exists B \, C,\eta B\,\coq{balls} \, \wedge \eta \, C\,\coq{balls} \, \wedge \\
&\eta \, P\, (\GBasicPointSetName \, B) \, \wedge \, \eta \, Q\, (\GBasicPointSetName \, C) \wedge \, \eta \, B\, (\coq{ext} C).
\end{aligned}
\]
\end{lemma}
The last lemma asserts that distinct balls generate corresponding points plurals in which some of them are disjoints:

\begin{lemma}[\coq{distinct\_centers\_separated\_by\_basic\_balls}]
\label{extbaseballs}
\[
\begin{aligned}
\forall P \, Q,  \eta \, P\,\coq{balls}\ \rightarrow\  &\eta \, Q\,\coq{balls}\ \rightarrow\ \neg\SameCenter\, P \, Q \rightarrow\ \exists B \, C,\eta B\,\coq{balls} \, \wedge \eta \, C\,\coq{balls} \, \wedge \\
&\eta \, P\, (\GBasicPointSetName \, B) \, \wedge \, \eta \, Q\, (\GBasicPointSetName \, C) \, \wedge \\ &\neg\PointSetMeet (\GBasicPointSetName B) (\GBasicPointSetName C).
\end{aligned}
\]
\end{lemma}

The formal sequence
\[
\begin{aligned}
&\coq{exterior\_GBasicPointSet\_disjoint},\\
&\coq{distinct\_centers\_have\_exterior\_basic\_balls},\\
&\coq{distinct\_centers\_separated\_by\_basic\_balls}.
\end{aligned}
\]
therefore gives:
\begin{theorem}[\coq{GPointHausdorff\_from\_Tarski\_separation}]
\label{thm:gpointhausdorff}
\[
P\not\sim Q\ \longrightarrow\
\exists u\,v,\ \GPointOpenName(u)\wedge\GPointOpenName(v)\wedge
\eta P\,u\wedge\eta Q\,v\wedge\neg\PointSetMeet(u,v).
\]
\end{theorem}
Thus \((X_G,\tau_G)\) is Hausdorff. This use of
\coq{Three\_points} is separate from the conditional non-discreteness
result.

\subsection{Regional-to-point comparison}
\label{subsec:phi-on-balls}

The earlier regional-to-point translation is
\[
\Phi(Q)=\PointClassInteriorName(\coq{pt}\,Q).
\]
For a ball generator, its image agrees exactly with the corresponding
basic geometric plural:
\begin{theorem}[\coq{Phi\_on\_balls}]\label{thm:phi-on-balls}
\[
\forall Q,\quad
\eta Q\,\coq{balls}\ \longrightarrow\
\PointSetEq(\GBasicPointSetName(Q),\Phi(Q)).
\]
\end{theorem}
The proof unfolds \(\Phi\) and applies
\coq{GBasicPointSet\_as\_interior}; the premise in the definition of
\PointSetEq{} already supplies a ball representative. Consequently, the
geometric basis is the point-class image of the family of mereological
balls:
\[
\{\,\GBasicPointSetName(Q)\mid\eta Q\,\coq{balls}\,\}
 =
\{\,\Phi(Q)\mid\eta Q\,\coq{balls}\,\}.
\]
This is a comparison of generators. It neither claims that \(\Phi\)
preserves arbitrary regional joins nor asserts that \(\Phi\) is a frame
morphism.

\section{Geometric Adherence and Kuratowski Closure}
\label{sec:geometric-closure}

\subsection{Geometric neighbourhoods}
\label{subsec:geometric-neighbourhoods}

The point-class topology defines a neighbourhood relation by
\begin{definition}[\GPointNeighbourhoodName]\label{def:gneighbourhood}
\[
\begin{aligned}
\GPointNeighbourhoodName(P,v)
\quad\triangleq\quad&
\eta P\,\coq{balls}\wedge\exists u,\\
&\GPointOpenName(u)\wedge\eta P\,u \,\wedge \, \PointSetIncl(u,v).
\end{aligned}
\]
\end{definition}
Whenever \(P\in\GBasicPointSetName(Q)\), the basic plural itself is a geometric neighbourhood \(v\) of \(P\); this result appears in the lemma \coq{GBasicPointSet\_is\_neighbourhood}: 
\begin{lemma}[\coq{GBasicPointSet\_is\_neighbourhood}]
\label{baseneighbour}
\[
\begin{aligned}
\forall P \, Q,  \eta \, P\, \GBasicPointSetName(Q) \rightarrow\  \GPointNeighbourhoodName (P \, (\GBasicPointSetName Q)).
\end{aligned}
\]
\end{lemma}

More generally,
\coq{GPointOpen\_neighbourhood} provides a geometric open neighbourhood
for every ball representative.

\subsection{Geometric closure}
\label{subsec:geometric-gclosure}

For a point-class plural \(u\), define \(\GClosurePointName(u)\) by the
usual neighbourhood condition:
\begin{definition}[\GClosurePointName]\label{def:gclosure}
\[
\eta P\,\GClosurePointName(u)
\quad\triangleq\quad
\eta P\,\coq{balls}\wedge
\forall v,\ \GPointOpenName(v)\rightarrow\eta P\,v
\rightarrow\PointSetMeet(v,u).
\]
\end{definition}
Here \PointSetMeet{} is extensional: two plurals meet when they contain
representatives of the same concentric point-class. Thus the definition
expresses adherence of reconstructed locations, rather than syntactic
identity of individual representatives.

The closure construction is well defined at the point-class level:
\begin{lemma}[\coq{GClosurePoint\_stable}]
\label{lem:gclosure-stable}
\[
\forall u,\quad \PointSetStable(\GClosurePointName(u)).
\]
\end{lemma}
\begin{lemma}[\coq{GClosurePoint\_respects\_PointSet\_eq}]
\label{lem:gclosure-respects-pointseteq}
\[
\PointSetEq(u,w)\ \longrightarrow\
\PointSetEq(\GClosurePointName(u),\GClosurePointName(w)).
\]
\end{lemma}
Thus closure maps arbitrary plurals to stable point-class plurals and
depends only on their extensional point-set content.

\subsection{The four Kuratowski axioms}
\label{subsec:gkuratowski}

\begin{theorem}[\coq{GClosurePoint\_is\_Kuratowski}]
\label{thm:gclosure-kuratowski}
\[
\begin{array}{ll}
\mathrm{(K0)}&
\GClosurePointName(\Lambda)=_{pt}\Lambda,\\[1mm]
\mathrm{(K1)}&
u\subseteq_{pt}\GClosurePointName(u),\\[1mm]
\mathrm{(K2)}&
\GClosurePointName(u\cup_{pt}v)
=_{pt}\GClosurePointName(u)\cup_{pt}\GClosurePointName(v),\\[1mm]
\mathrm{(K3)}&
\GClosurePointName(\GClosurePointName(u))
=_{pt}\GClosurePointName(u).
\end{array}
\]
\end{theorem}
The four clauses are proved by
\coq{GClosurePoint\_empty}, \coq{GClosurePoint\_extensive},
\coq{GClosurePoint\_union}, and \coq{GClosurePoint\_idempotent}.
For the union law, two open witnesses disjoint from \(u\) and \(v\) can be
intersected to obtain a witness disjoint from their union. Idempotence
uses the stability of geometric opens under replacement of a representative
by a concentric representative. These proofs need neither
\coq{Three\_points} nor \coq{BallCenterBundle}; those hypotheses concern
the additional separation and non-discreteness properties above.

\subsection{Limitations and Representation Problem}
\label{subsec:gclosure-scope}

Since Theorem~\ref{thm:gpointopen-topology} is established independently,
\GClosurePointName{} is not merely an abstract Kuratowski operator. It is
the neighbourhood closure induced by the ball-generated geometric topology
\(\tau_G\). The paper therefore establishes the progression
\[
\begin{gathered}
\text{balls}\ \longrightarrow\ \text{basic point-plurals}
\ \longrightarrow\ \tau_G,\\
\tau_G\ \longrightarrow\ T_2\text{ and conditional non-discreteness}
\ \longrightarrow\ \GClosurePointName.
\end{gathered}
\]
Boundary, cavity, and hole constructions are deliberately left for later
work: the present result is the verified construction of the topology, its
separation properties, and its induced closure.

The construction established here is internal to the nominal Tarskian
framework. It proves that the quotient setoid
\(X_G=\mathrm{balls}/{\sim}\) carries a Hausdorff topology generated by
the plurals \coq{GBasicPointSet(Q)}, and that \GClosurePointName{} is the
corresponding Kuratowski closure operator. It does not yet prove that this
topology is the Euclidean topology in the intended standard model.

The natural standard target of the present construction is
three-dimensional Euclidean space. A Pieri structure is a pair
\(\langle\mathcal P,\triangle\rangle\), where \(\mathcal P\) is a
non-empty domain of abstract points and
\(\triangle\subseteq\mathcal P^{3}\) is the ternary equidistance
relation \cite{Pieri1908,GruszczynskiPietruszczak2007}. The formula
\[
x\,y\mathrel{\triangle}z
\]
means that \(x\) and \(y\) are equidistant from \(z\). In the standard
Euclidean interpretation this relation is defined by
\[
x\,y\mathrel{\triangle_{\mathbb R^{3}}}z
\quad\Longleftrightarrow\quad
d(x,z)=d(y,z),
\]
where \(d\) is the usual Euclidean distance.

Pieri's categoricity theorem states that every complete Pieri structure
is isomorphic to
\(\langle\mathbb R^{3},\triangle_{\mathbb R^{3}}\rangle\)
\cite{Pieri1908,GruszczynskiPietruszczak2007}.
Consequently, if \coq{Equid} induces a well-defined relation
\(\widehat{\mathrm{Equid}}\) on
\(X_G=\mathrm{balls}/{\sim}\), and if
\[
\langle X_G,\widehat{\mathrm{Equid}}\rangle
\]
is proved to satisfy the complete Pieri axiom system, the theorem supplies
a geometric isomorphism
\[
f:X_G\longrightarrow\mathbb R^{3}.
\]

This premise is not yet established in the present Coq development.
Although the library defines \coq{Equid}, proves several of its properties,
and assumes \coq{Three\_points}, these results do not by themselves show
that the quotient structure satisfies all the axioms of Pieri's geometry,
in particular its continuity and dimension axioms. The gap is therefore one
of axiomatic strength rather than merely a remaining translation task: the
current Coq signature contains neither the required continuity principles
nor axioms fixing the dimension to three.

A further compatibility theorem would then be required between the
geometric isomorphism \(f\) and the ball-generated topology. Let
\(\lbrack\!\lbrack Q\rbrack\!\rbrack\) denote the Euclidean open ball
interpreting the nominal ball representative \(Q\). The expected
basis-level correspondence is
\[
f\!\left(
 \left\{[B]_{\sim}\in X_G\mid
        \eta B\,\coq{GBasicPointSet}(Q)\right\}
\right)
=
\lbrack\!\lbrack Q\rbrack\!\rbrack
=
\operatorname{Int}_{\mathbb R^{3}}
  (\lbrack\!\lbrack Q\rbrack\!\rbrack).
\]
Here \(\operatorname{Int}_{\mathbb R^{3}}\) denotes topological interior
with respect to the standard Euclidean topology. Since
\(\lbrack\!\lbrack Q\rbrack\!\rbrack\) is an open Euclidean ball, its
interior is the ball itself. The displayed equality states that the
reconstructed points represented by
\coq{GBasicPointSet(Q)} correspond exactly to the Euclidean points lying
in the ball represented by \(Q\).

The occurrence of \(f\) on the left-hand side is essential: before its
application, the displayed set is a subset of the quotient \(X_G\);
after its application, it becomes a subset of \(\mathbb R^{3}\), and can
therefore be compared with
\(\lbrack\!\lbrack Q\rbrack\!\rbrack\).

A related representation result is Bennett's axiomatisation of
Region-Based Geometry (RBG). The notation \(RBG_n\) denotes the
\(n\)-dimensional member of this family of region-based theories; the
relevant specialization here is therefore \(RBG_3\). Bennett's
categoricity theorem associates reconstructed points with Cartesian
points, sphere representatives with Euclidean open balls, and regions
with non-empty regular open subsets of \(\mathbb R^{3}\). It provides a
second possible route to a representation theorem, but it cannot yet be
applied to the present development because no interpretation of the
nominal setoid structure into \(RBG_3\) has been formally established.

\section{Conclusion and Future Work}

The originality of the present work does not lie in quotient setoids,
basis-generated topologies, or Kuratowski closure taken separately. It
lies in their combination within Tarski's reconstruction of points and a
Le{\'s}niewskian nominal mereology. To the best of our knowledge, no
previous formal development has shown that concentric families of ball
representatives support a Hausdorff topology and a neighbourhood-based
Kuratowski closure while reconstructed points remain unreified in the
object ontology.

The construction identifies the precise quotient structure required for
this purpose: ball representatives are identified by \texttt{same\_center},
point-plurals are required to be stable under replacement by concentric
representatives, and regional balls generate a basis of such plurals. The
resulting topology is defined metatheoretically on the setoid, is Hausdorff
under Tarski's separation principles, and is non-discrete under an explicit
local richness condition. Its closure operator satisfies the four
Kuratowski axioms, maps arbitrary plurals to stable point-class plurals,
and respects extensional point-set equality.

This differs from classical reconstructions of Tarski's geometry conducted
in a set-theoretic setting, from contact-algebraic approaches that begin
with an algebra of regions, and from localic approaches that take the
lattice of opens as primitive. The contribution is therefore both
conceptual and formal: it isolates the quotient and invariance conditions
that make Tarski's reconstructed points topologically usable, and it
verifies the resulting construction in Coq without adding point individuals
or leaving proof obligations admitted. This result is internal to the
nominal framework; its identification with the Euclidean topology in the
standard model remains a separate representation problem.

We have developed a ball-generated topology on the reconstructed points of
Tarski's geometry of solids within the \(\lambda\)-MM nominal framework.
The construction separates regional ball generators, their concentric
point-classes, and stable plurals of their representatives.  It proves that
the resulting family of geometric opens is a topology on the setoid
\(\coq{balls}/{\sim}\).

The main result is a verified progression from the geometry of balls to
topology and closure.  The basic plurals cover the reconstructed-point
domain and refine local intersections; the resulting topology is Hausdorff
under the separation axiom \coq{Three\_points} and non-discrete under the explicit
\coq{BallCenterBundle} condition.  Its neighbourhood closure
\GClosurePointName{} satisfies the four Kuratowski axioms.  The theorem
\coq{Phi\_on\_balls} identifies the basic geometric plurals with the
point-class images of regional ball generators.

The formalization does not introduce a second ontology.  Both the regional
and point-class presentations are defined over the same nominal domain of
\(\lambda\)-MM.  The passage to point-classes is a change of topological
representation imposed by the reconstruction of points from concentric
balls, not an ontological duplication.

\paragraph{Future work.}
The next step is to develop boundary, cavity, and hole constructions from
the geometric topology just established, and then to investigate their
topological semantics for spatial regions.  Rather than viewing a region
solely as an ontological individual or mereological whole, one may
associate with it a derived structure through its induced point-class
plurals.  This direction is related to the formal-ontological analysis of
spatial regions and boundaries developed by Smith and Varzi
\cite{SmithVarzi2000,Smith1995}. Their work provides a rich classification
of boundaries; the present framework suggests a complementary programme in
which neighbourhood, closure, and boundary operations receive a formally
verified topological semantics.

This programme would bridge nominal mereology, point-free geometry, and
formal ontology while preserving a common nominal domain.

\end{document}